\documentclass[reqno]{amsart}
\usepackage[english]{babel}
\usepackage{amscd,amssymb,amsmath,amsfonts,latexsym,amsthm}
\usepackage{inputenc}
\usepackage{graphicx}
\usepackage{hyperref,amsthm,amsmath,amssymb,color,multirow,graphicx}
\usepackage[fontsize=11pt]{scrextend}

\newtheorem{thm}{Theorem}[section]
\newtheorem{cor}[thm]{Corollary}
\newtheorem{lem}[thm]{Lemma}
\newtheorem{prop}[thm]{Proposition}
\theoremstyle{definition}
\newtheorem{defn}[thm]{Definition}
\newtheorem{rem}[thm]{Remark}
\newtheorem{ex}[thm]{Example}

\numberwithin{equation}{section}

\def\e{{\bf 1}\!\!{\rm I}}

\newcommand{\be}{\begin{equation}}

\newcommand{\ee}{\end{equation}}

\begin{document}

\begin{center}
\textbf{\large  Algebras of $m$-measurable and locally $m$-measurable operators affiliated with a W*-algebra and their derivations. }\\

\textsc{\textbf{Rakhimov Abdugafur Abdumajidovich}}\\
\textsc{National University of Uzbekistan named after M.Ulugbek, Tashkent, Uzbekistan\\
Karadeniz Technical University, Trabzon, Turkey.}\\
{rakhimov@ktu.edu.tr}\\
\textsc{\textbf{Tukhtasinov Tokhirjon Shokirjon o'g'li }}\\
\textsc{Fergana State University, Fergana, Uzbekistan}\\
{toxirjontoxtasinov13@gmail.com}\\

\end{center}



\makeatletter
\renewcommand{\@evenhead}{\vbox{\thepage \hfil {\it Rakhimov A.A., Tukhtasinov T.Sh.}  \hrule }}
\renewcommand{\@oddhead}{\vbox{\hfill
{\it  Algebras of $m$-measurable and locally $m$-measurable operators affiliated with a W*-algebra. }\hfill
\thepage \hrule}} \makeatother

\label{firstpage}

\section{Introduction}

In \cite{Ciach}, Leszek J. Ciach introduced the concept of
$m$-measurable operators affiliated with von Neumann
algebras and studied their properties.

In this paper, we investigate the relationships between
the algebras of $m$-measurable and locally $m$-measurable
operators and the classical algebras of measurable and
locally measurable operators.

Furthermore, we study derivations on these algebras,
with particular attention to the question of whether
they are inner.

\section{Preliminaries}

Let $M$ be an algebra and let $P(M)$ be the set of all projections of $M$.
Recall that projectors $p$ and $q$ are called {\it equivalent} if there exists a partial isometry $v\in M$ such that
$v^*v=p$ and $vv^*=q$ and is denoted by $p\sim q$; say that $p\leq q$, if $pq=p$.

By {\it subadditive measure} on $M$ we mean a mapping $m:P(M)\to [0,\infty]$ with the following properties (see \cite{Ciach}):

i) $m(0)=0$, and $m(p)=0$ implies $p=0$ (faithfulness);

ii) $p\le q$ implies $m(p)\le m(q)$ (monotonicity);

iii) $p\sim q$ implies $m(p)=m(q)$;

iv) $m(p\vee q) \le m(p) + m(q)$ (subadditivity);

v) $p_n\uparrow p$ implies $m(p_n)\uparrow m(p)$.

\def\e{{\bf 1}\!\!{\rm I}}

Now, let $B(H)$ be the algebra of all bounded linear operators on a complex Hilbert space $H$. A weakly closed *-subalgebra
$M$ containing the identity operator $\e$ in $B(H)$ is called a W*-{\it algebra}.
The set $M'=\{x\in B(H): xy=yx, \forall y\in M\}$ is called the {\it commutant} of $M$.
A *-algebra $M$ is called the {\it von Neumann algebra}, if $M'' = M$.
It is known that a *-algebra $M$ is a von Neumann algebra if and only if it is a W*-algebra.
$Z = M \cap M'$ is called the \textit{center} of $M$.(\cite{ARU}, \cite{AR}, \cite{Li})\\[-2mm]

\begin{defn} 1) An operator $a$ acting in a Hilbert space $H$ with a domain of definition $Dom(a)$ is called {\it affiliated with the algebra} $M$ (symbolically, $a\eta M$) if $au'=u'a$, for any unitary operator $u'\in M'$.

2) A linear subspace $D\subset H$ is called {\it affiliated with the algebra} $M$ if $u'(D)\subset D$ for any unitary operator $u'\in M'$ and denoted by $D\eta M$.

3) A linear subspace $D\subset H$ with $D\eta M$ is said to be {\it $m$-dense} in $H$ if for any $\varepsilon >0$ there exists a projection $p\in P(M)$ such that $p(H)\subset D$ and $m(p^\perp)\leq \varepsilon$.

4) An operator $a$ is called {\it $m$-measurable}, if (i) $a\eta M$; (ii) $D(a)$ is $m$-dense; (iii) $a$ is closed.

The set of all $m$-measurable operators is denoted by $L_m(M)$.

5) An operator $a$ is called {\it locally $m$-measurable} if $a\eta M$ and there exists a sequence of central projections
$\{z_n\}_{n=1}^{\infty}$ such that $z_n\uparrow \e$ and, for every $n\in\mathbb N$, $z_na\in L_m(M)$.

The set of all locally $m$-measurable operators is denoted by $LL_m(M)$.

6) A closed densely defined operator $x$ is called \emph{measurable with respect to $M$} if \ $x\eta M$
and its domain $D(x)$ is strongly dense in $H$ with respect to $M$.
The set of all measurable operators with respect to $M$ is denoted by
$S(M)$.

7) A closed densely defined operator $x$ affiliated with $M$ is called
\emph{locally measurable with respect to $M$} if there exists a sequence of central projections
\ $\{z_n\}_{n=1}^{\infty}$ such that \ $z_n\uparrow \e$ and $z_nx\in S(M)$ for $n\in\mathbb N$.
The set of all locally measurable operators affiliated with $M$ is denoted by $LS(M)$ (\cite{Segal}, \cite{MCH}, \cite{AlAyK}).
\end{defn}

AI disclosure. Errors in a previous version were found by GPT 6 Astra.

\section{Relationship between algebras $L_m(M)$, $LL_m(M)$ and $S(M)$, $LS(M)$.}

We will show that, in the general case, the following holds $L_m(M)\not=S(M)$ and $LL_m(M)\not=LS(M)$.
Let us demonstrate this through several examples.\\[-2mm]

\begin{ex}
Consider the von Neumann algebras $ M_1=L^\infty(0,1)\,\bar\otimes\,B(\ell^2)$ and $M_2=L^\infty(0,1)$. Let $M=M_1\oplus M_2$.
Then $M$ is a type I von Neumann algebra, and every projection $p\in P(M)$ has the form $p=p_1\oplus p_2$, where $p_1\in P(M_1)$ and $p_2\in P(M_2)$. Let $c(p_1)$ denote the central support of $p_1$ in $M_1$, and let $\mu$ be the Lebesgue measure on $(0,1)$.
Define a map $m:P(M)\longrightarrow[0,\infty]$ by
$$
m(p_1\oplus p_2)=\mu(c(p_1))+\mathbf{\e}_{\{p_2\neq0\}},
$$
where
$$
\mathbf{\e}_{\{p_2\neq0\}}=
\begin{cases}
0, & p_2=0,\\[1mm]
1, & p_2\neq0.
\end{cases}
$$
Equivalently, define $m_1(p_1)=\mu(c(p_1))$ and
$$
m_2(p_2)=
\begin{cases}
0, & p_2=0,\\[1mm]
1, & p_2\neq0.
\end{cases}
$$
Then
$$
m(p_1\oplus p_2)=m_1(p_1)+m_2(p_2).
$$

It is straightforward to verify that $m$ is a faithful, finite, subadditive measure on the
projection lattice $P(M)$.\\[-2mm]

Consider the affiliated operator $x=\left(\frac{1}{t}I_{\ell^2}\right)\oplus 0$, where $I_{\ell^2}=id_{\ell^2}$. Then we have $x\in L_m(M)$ and $x\notin S(M)$, hence $L_m(M)\not\subseteq S(M)$.\\[-2mm]

Now let us show the converse. We consider the affiliated operator $y=0\oplus\frac{1}{t}$.\\[-2mm]

Since $M_2=L^\infty(0,1)$ is a finite von Neumann algebra, the operator $t\mapsto1/t$ belongs to $S(M_2)$.
Moreover, measurable operators over a finite direct sum of von Neumann algebras are characterized componentwise.
Therefore, $y\in S(M)$. But we have $y\notin L_m(M)$, and hence $S(M)\not\subseteq L_m(M)$.\end{ex}
\begin{ex} Let $M$ be a type II$_1$ factor.
Define a mapping $m:P(M)\longrightarrow[0,\infty]$ by
$$m(p)=
\begin{cases}
0, & p=0,\\
1, & p\neq0.
\end{cases}
$$
It is clear that $m$ is a subadditive measure. Then $LL_m(M)=M\subsetneq LS(M)$. Thus $LS(M)\not\subseteq LL_m(M)$.\end{ex}

%\begin{center}
%{\Large\bfseries Stage I: Characterization of Subadditive Measures\\ on Projections of Type $\mathrm{II}_1$ Factors}\par\medskip
%{\small A complete trace-parametrization theorem and its consequences}
%\end{center}

\section{Characterization of subadditive measures on projections of type II$_1$ Factors.}

\def\e{{\bf 1}\!\!{\rm I}}

%\section*{Setup}
Let $R$ be a type II$_1$ factor, with faithful normal
normalized trace $\tau$ ($\tau(\e)=1$). No separability assumption is
made. Denote by $P(R)$ the lattice of projections of $R$.

%The standard dimension theory of a type $\mathrm{II}_1$ factor yields
%\begin{equation}\label{eq:dimension}
%\tau(P(R))=[0,1],\qquad
%p\sim_{R}q\ \Longleftrightarrow\ \tau(p)=\tau(q).
%\end{equation}
%See, for example, Anantharaman--Popa~\cite[Proposition~4.1.6]{AP}. We begin with a useful strengthened form of the range assertion.

\begin{lem}[An increasing trace scale]\label{lem:trace-scale}
There exists a family $(p_t)_{0\leq t\leq1}\subseteq P(R)$ such that
$\tau(p_t)=t$, \ $p_s\leq p_t$ \ ($s\leq t$),
and, for every sequence $t_n\uparrow t$, we have $p_{t_n}\uparrow p_t$.
\end{lem}
\begin{proof}
Every nonzero corner of a type $\mathrm{II}_1$ factor is diffuse.
By the halving property, each projection can be split into two mutually
orthogonal equivalent subprojections. Iterating, construct coherent
dyadic partitions
\[
\e =\sum_{k=1}^{2^n}e_{n,k},\qquad
\tau(e_{n,k})=2^{-n},\qquad
 e_{n,k}=e_{n+1,2k-1}+e_{n+1,2k}
\]
for $n\geq0$, with $e_{0,1}=1$.
For a dyadic rational $r=k2^{-n}\in[0,1]$, set
\[
p_r=\sum_{j=1}^{k}e_{n,j},\qquad p_0=0.
\]
Coherence makes this definition independent of the dyadic
representation. These projections are increasing in $r$, and
$\tau(p_r)=r$.
For arbitrary $t\in(0,1]$ let
\[
p_t=\bigvee_{\substack{r\in\mathbb D\\0\le r<t}}p_r,
\qquad
\mathbb D=\left\{k2^{-n}:n\geq0,\ 0\le k\le2^n\right\}.
\]
Normality of $\tau$ gives $\tau(p_t)=t$. If $t_n\uparrow t$,
then every dyadic $r<t$ satisfies $r<t_n$ for all sufficiently
large $n$. Thus $\bigvee_n p_{t_n}=p_t$, as required.
\end{proof}

\begin{thm}[Complete characterization of subadditive measures]
\label{thm:characterization}
Let $R$ be an arbitrary type $\mathrm{II}_1$ factor with normalized
trace $\tau$. The following statements are equivalent.
\begin{enumerate}
\item[(a)] $m:P(R)\to[0,\infty]$ is a subadditive measure.
\item[(b)] There exists a \emph{unique} function
$f:[0,1]\to[0,\infty]$ such that
\begin{equation}\label{eq:factorization}
m(p)=f(\tau(p))\qquad(p\in P(R)),
\end{equation}
and the following four conditions hold:
\begin{align}
f(0)&=0,\qquad f(t)>0\quad(0<t\leq1); \label{eq:f1}\\
f(s)&\leq f(t)\quad(0\leq s\leq t\leq1);\label{eq:f2}\\
f(s+t)&\leq f(s)+f(t)
\quad(s,t\geq0,\ s+t\leq1);\label{eq:f3}\\
f(t_n)&\uparrow f(t)\quad\text{whenever }t_n\uparrow t
\text{ in }[0,1].\label{eq:f4}
\end{align}
\end{enumerate}
Consequently, the assignment $m\leftrightarrow f$ is a bijection
between Ciach subadditive measures on $P(R)$ and functions
satisfying \eqref{eq:f1}--\eqref{eq:f4}. Condition
\eqref{eq:f4} is equivalent to left-continuity of $f$ at
every $t\in(0,1]$, in the extended-real sense.
\end{thm}
\begin{proof}
\emph{(a) $\Rightarrow$ (b).}
Let $(p_t)$ be the trace scale of Lemma \ref{lem:trace-scale}, and define
$f(t)=m(p_t)$, \ $0\leq t\leq 1$.
For every $p\in P(R)$, equality of traces gives
$p\sim_{R}p_{\tau(p)}$. Therefore, $m(p)=m(p_{\tau(p)})=f(\tau(p))$.
Since every $t\in[0,1]$ is the trace of some projection,
$f$ is uniquely determined by $m$.

Faithfulness proves \eqref{eq:f1}, and the monotonicity of the
trace scale gives \eqref{eq:f2}. For $s,t\geq0$ with $s+t\leq1$,
the projections $p_s$ and $q=p_{s+t}-p_s$
are orthogonal, with $\tau(q)=t$. Hence $q\sim_{R}p_t$, and
subadditivity together with invariance yields
\[
\begin{aligned}
f(s+t)&=m(p_{s+t})=m(p_s\vee q)\\
&\leq m(p_s)+m(q)=f(s)+f(t).
\end{aligned}
\]
Thus \eqref{eq:f3} holds. Finally, if $t_n\uparrow t$, then
$p_{t_n}\uparrow p_t$ by Lemma \ref{lem:trace-scale}.
By continuity-from-below axiom, $f(t_n)=m(p_{t_n})\uparrow m(p_t)=f(t)$, which is \eqref{eq:f4}.

\medskip\noindent
\emph{(b) $\Rightarrow$ (a).}
Suppose that $f$ satisfies \eqref{eq:f1}--\eqref{eq:f4}, and set
$m(p)=f(\tau(p))$. Faithfulness follows from \eqref{eq:f1}
and faithfulness of $\tau$. Monotonicity follows from
\eqref{eq:f2}. Murray--von Neumann invariance follows
from the tracial property of $\tau$.

To prove subadditivity for \emph{arbitrary} $p,q\in P(R)$,
recall the dimension identity
\begin{equation}\label{eq:modular}
\tau(p\vee q)+\tau(p\wedge q)=\tau(p)+\tau(q).
\end{equation}
Indeed, $(p\vee q)-p\sim_{R}q-(p\wedge q)$, and taking
traces gives \eqref{eq:modular}.
Put
\[
s=\tau(p),\qquad
r=\tau(p\vee q)-\tau(p)
=\tau(q)-\tau(p\wedge q).
\]
Then $0\leq r\leq\tau(q)$ and $s+r=\tau(p\vee q)\leq1$.
By \eqref{eq:f2} and \eqref{eq:f3},
\[
\begin{aligned}
m(p\vee q)&=f(s+r)\\
&\leq f(s)+f(r)\\
&\leq f(\tau(p))+f(\tau(q))
=m(p)+m(q).
\end{aligned}
\]
Finally, if $p_n\uparrow p$, normality of the trace yields
$\tau(p_n)\uparrow\tau(p)$. By \eqref{eq:f4},
\[
m(p_n)=f(\tau(p_n))\uparrow f(\tau(p))=m(p).
\]
Thus, $m$ satisfies all the axioms of subadditivity.
\end{proof}

\begin{cor}[The small-projection parameter]
\label{cor:delta}
For a subadditive measure $m$ associated with $f$, define
\[
\delta_m=\inf\{m(p):p\in P(R),\ p\neq0\}.
\]
Then
\begin{equation}\label{eq:delta}
\delta_m=\inf_{0<t\leq1}f(t) = \lim_{t\downarrow0}f(t)\in[0,\infty].
\end{equation}
The axioms of a subadditive measure do \emph{not} force
$\delta_m=0$.
\end{cor}
\begin{proof}
Every $t\in(0,1]$ is the trace of a nonzero projection, and
every nonzero projection has strictly positive trace. Hence the two
infima agree. The limit equality follows from monotonicity of $f$.
The function $f(t)=\e_{\{t>0\}}$ satisfies all four
conditions of Theorem \ref{thm:characterization} but has
$\delta_m=1$.
\end{proof}

\begin{cor}[Finiteness and continuity]\label{cor:finite}
For every subadditive measure $m$ on a type $\mathrm{II}_1$ factor,
exactly one of the following holds:
\begin{enumerate}
\item $m(p)=\infty$ for every nonzero $p\in P(R)$;
\item $m(\e)<\infty$, so $f$ is finite-valued and bounded on $[0,1]$.
\end{enumerate}
If $\delta_m=0$, the second alternative necessarily holds and
$f$ is continuous on $[0,1]$.
\end{cor}
\begin{proof}
If $f(t_0)<\infty$ for some $t_0>0$, choose an integer $N$
with $1/N\leq t_0$. Monotonicity and repeated use of
\eqref{eq:f3} yield
\[
f(1)\leq Nf(1/N)\leq Nf(t_0)<\infty.
\]
Thus, unless $f$ is identically infinite on $(0,1]$, it is
finite-valued and bounded by $f(1)$.

If $\delta_m=0$, then $f(h)\to0$ as $h\downarrow0$,
so the finite-valued alternative applies. For $0\leq t<t+h\leq1$,
monotonicity and subadditivity give
\[
0\leq f(t+h)-f(t)\leq f(h)\xrightarrow[h\downarrow0]{}0.
\]
Thus $f$ is right-continuous on $[0,1)$ and at $0$.
Left-continuity at each $t\in(0,1]$ follows from
\eqref{eq:f4}. Hence $f$ is continuous on $[0,1]$.
\end{proof}

\begin{ex}[Trace measure]
$f(t)=t$ corresponds to $m(p)=\tau(p)$, and $\delta_m=0$.
\end{ex}
\begin{ex}[Nonlinear trace measures]
For $0<\alpha\leq1$, the function $f(t)=t^\alpha$ satisfies
\eqref{eq:f1}--\eqref{eq:f4} and yields
$m(p)=\tau(p)^\alpha$. Here $\delta_m=0$.
\end{ex}
\begin{ex}[Indicator measure]
Define $f(0)=0$ and $f(t)=1$ for $t>0$. Then
\[
m_0(p)=\e_{\{p\neq0\}},\qquad\delta_{m_0}=1.
\]
This example is continuous from below on projections, although
$f$ has a jump at $0$. Thus Ciach's normality axiom is not
equivalent to continuity from above at zero.
\end{ex}
\begin{ex}[Mixed and infinite-valued measures]
For $0<\alpha\leq1$, the function
$f(t)=\e_{\{t>0\}}+t^\alpha$ satisfies the conditions,
and $\delta_m=1$.
The function $f(0)=0$, $f(t)=\infty$ for $t>0$ also satisfies
the conditions, with $\delta_m=\infty$.
\end{ex}

\begin{rem}[Scope and next step]
The representation $m=f\circ\tau$ rests on the classical
classification of projections by trace in type $\mathrm{II}_1$
factors. It is a structural starting point, not by itself a claim
of a new derivation theorem. In Stage II we will use the parameter
$\delta_m$ to determine precisely when $L_m(R)=R$ and when
$L_m(R)=S(R)$, and then examine $LL_m(R)$.
\end{rem}

%\begin{center}
%{\Large\bfseries Stage II: A Dichotomy for Locally\
%$m$-Measurable Operators on Type $\mathrm{II}_1$ Factors}\par\medskip
%{\small A complete classification, its converse, and explicit separating examples}
%\end{center}

\section{Locally $m$-measurable operators on type II$_1$ factors.}

%\section*{1. Definitions and results from Stage I}

Let $R\subseteq B(H)$ be a type II$_1$ factor,
not necessarily separable, and let $\tau$ be its faithful normal
tracial state, normalized by $\tau(\e)=1$.
Let $m: P(R)\to[0,\infty]$ be a subadditive measure.
In particular, $m$ is faithful, monotone, invariant under
Murray--von Neumann equivalence, subadditive on projection joins,
and continuous from below.

Write $Aff(R)$ for the set of all closed densely defined operators
affiliated with $R$; write $S(R)$ and $LS(R)$ for the Segal-measurable
and locally measurable operators, respectively. All products of
unbounded affiliated operators are taken with the usual closed-operator
interpretation.
\begin{lem}\label{lem:finite-factor}
For every type $\mathrm{II}_1$ factor $R$ one has
\[
S(R)=LS(R)=Aff(R),\qquad LL_m(R)=L_m(R).
\]
Moreover, $R\subsetneq S(R)$.
\end{lem}

\begin{proof}
As $R$ is finite, every projection of $R$ is finite in the
Murray--von Neumann sense. Let $x\in Aff(R)$, and set
$p_n=E_{[0,n]}(|x|)$. By the spectral theorem, $p_n\uparrow\e$,
$p_nH\subseteq\mathcal D(x)$, and $\e-p_n$ is finite.
Consequently, $\mathcal D(x)$ is strongly dense and $x\in S(R)$.
The converse $S(R)\subseteq Aff(R)$ holds by definition, proving
$S(R)=Aff(R)$.

Because $R$ is a factor, $ P(Z(R))=\{0,1\}$. If $(z_n)$ is an
increasing sequence of central projections converging to $\e$, then
$z_n=\e$ for all sufficiently large $n$. Thus $LS(R)=S(R)$ and $LL_m(R)=L_m(R)$.

To see that $R\neq S(R)$, recursively halve the identity:
put $r_0=1$ and choose orthogonal nonzero projections $e_n,r_n\in R$
such that
\[
r_{n-1}=e_n+r_n,\qquad
\tau(e_n)=\tau(r_n)=2^{-n},\qquad n\ge1.
\]
The projections $e_n$ are mutually orthogonal and
$\tau(\bigwedge_n r_n)=\lim_n2^{-n}=0$. Faithfulness of $\tau$
implies $\bigwedge_n r_n=0$; hence $\bigvee_ne_n=\e$.
Define the positive operator
\begin{equation}\label{eq:witness}
h=\sum_{n=1}^{\infty}n e_n,
\qquad
\mathcal D(h)=\left\{\xi\in H:
\sum_{n=1}^{\infty}n^2\|e_n\xi\|^2<\infty\right\}.
\end{equation}
Its domain is dense, $h$ is self-adjoint and affiliated with $R$,
and it is unbounded because $he_n=ne_n$ and $e_n\neq0$ for all $n$.
Thus $h\in S(R)\setminus R$.
\end{proof}

\begin{lem}[Spectral criterion for $m$-measurability]
\label{lem:Ciach}
For every closed densely defined operator $x$ affiliated with $R$,
\begin{equation}\label{eq:spectral-criterion}
x\in L_m(R)\quad\Longleftrightarrow\quad
\lim_{t\to\infty}m\bigl(E_{(t,\infty)}(|x|)\bigr)=0.
\end{equation}
\end{lem}
\begin{proof}
This is precisely the combination of parts (i) and (ii) of
\cite[Proposition 2.1]{Ciach}. Part (i) reduces the assertion
to $|x|$, and part (ii) gives the stated spectral-tail criterion for
positive affiliated operators.
\end{proof}
%\section*{2. The dichotomy theorem}
\begin{thm}[Complete dichotomy]\label{thm:dichotomy}
Let $R$ be an arbitrary type $\mathrm{II}_1$ factor, and let
$m: P(R)\to[0,\infty]$ be a subadditive measure. Then, with no separability or finiteness assumption on $m$,
\begin{equation}\label{eq:main}
LL_m(R)=L_m(R)=
\begin{cases}
S(R)=Aff(R),&\delta_m=0,\\[2mm]
R,&\delta_m>0.
\end{cases}
\end{equation}
Here the second case includes $\delta_m=\infty$.
\end{thm}

\begin{proof}
Since $R$ is a factor, Lemma~\ref{lem:finite-factor} gives
\begin{equation}\label{eq:factor-trivial}
LL_m(R)=L_m(R),\qquad S(R)=Aff(R).
\end{equation}
It remains to identify $L_m(R)$.

\medskip\noindent
{\it Case 1}: $\delta_m=0$.
Take any $x\in S(R)=Aff(R)$, and write
\[
q_t=E_{(t,\infty)}(|x|),\qquad t>0.
\]
Since $|x|$ is self-adjoint and densely defined, the spectral theorem
gives $q_t\downarrow0$ strongly as $t\to\infty$. The trace $\tau$
is finite and normal; therefore
\begin{equation}\label{eq:trace-tail}
\tau(q_t)\longrightarrow0\qquad(t\to\infty).
\end{equation}
By \eqref{eq:factorization} and \eqref{eq:delta},
\[
m(q_t)=f(\tau(q_t))\longrightarrow0,
\]
because $f(s)\to\delta_m=0$ as $s\downarrow0$, including the
possibility that $\tau(q_t)=0$ at some finite stage (when $f(0)=0$).
By Lemma~\ref{lem:Ciach}, $x\in L_m(R)$. Thus
\[
S(R)\subseteq L_m(R).
\]
Conversely, by Ciach's definition, every element of $L_m(R)$ is
a closed densely defined operator affiliated with $R$, so
$L_m(R)\subseteq Aff(R)=S(R)$. Therefore $L_m(R)=S(R)\qquad(\delta_m=0)$.

\medskip\noindent
{\it Case 2}: $\delta_m>0$.
Let $x\in L_m(R)$ and again set $q_t=E_{(t,\infty)}(|x|)$.
By Lemma~\ref{lem:Ciach}, $m(q_t)\to0$.
Select a finite real number $\varepsilon$ with $0<\varepsilon<\delta_m$;
this is possible even if $\delta_m=\infty$.
There exists $t_0>0$ such that $m(q_{t_0})<\varepsilon$.
If $q_{t_0}\neq0$, the definition of $\delta_m$ gives
$m(q_{t_0})\ge\delta_m>\varepsilon$,
a contradiction. Hence $q_{t_0}=0$, so $E_{(t_0,\infty)}(|x|)=0$.
By the spectral theorem, $|x|\le t_0 1$, and $x$ is bounded.
Since it is affiliated with $R$, $x\in R$. Thus
$L_m(R)\subseteq R$.

Conversely, if $x\in R$, then $E_{(t,\infty)}(|x|)=0$ for all
$t\ge\|x\|$, so $x\in L_m(R)$ by Lemma~\ref{lem:Ciach}.
Therefore $L_m(R)=R$ \ $(\delta_m>0)$.
Combining both cases with \eqref{eq:factor-trivial} proves
\eqref{eq:main}.
\end{proof}

\begin{cor}[Necessary and sufficient conditions]
\label{cor:iff}
For a subadditive measure $m$ on a type $\mathrm{II}_1$ factor,
\[
\begin{aligned}
LL_m(R)=S(R)
&\quad\Longleftrightarrow\quad\delta_m=0,\\
LL_m(R)=R
&\quad\Longleftrightarrow\quad\delta_m>0.
\end{aligned}
\]
Furthermore, $LL_m(R)$ is never a proper intermediate algebra
between $R$ and $S(R)$.
\end{cor}
\begin{proof}
The forward implications follow from Theorem~\ref{thm:dichotomy}
and the strict inclusion $R\subsetneq S(R)$ proved in
Lemma~\ref{lem:finite-factor}. The reverse implications are exactly
Theorem~\ref{thm:dichotomy}.
\end{proof}

\begin{cor}[A universal separating affiliated operator]
\label{cor:witness}
Let $h$ be the unbounded positive operator in
\eqref{eq:witness}. Then
\[
E_{(N,\infty)}(h)=\bigvee_{n>N}e_n,
\qquad \tau(E_{(N,\infty)}(h))=2^{-N}
\quad(N\in\mathbb N),
\]
and
\begin{equation}\label{eq:tail-witness}
m(E_{(N,\infty)}(h))=f(2^{-N})\longrightarrow\delta_m.
\end{equation}
Consequently, \ $h\in L_m(R)\ \Longleftrightarrow\ \delta_m=0$.
\ The same operator $h$ therefore distinguishes the two classes of
subadditive measures.
\end{cor}
\begin{proof}
The spectral decomposition of $h$ and $\tau(e_n)=2^{-n}$ give the
stated formulas. By monotonicity of the spectral tails and
Lemma~\ref{lem:Ciach}, the last assertion follows.
\end{proof}
%\section*{3. Examples and implications for the next stages}

\begin{ex}[The normalized trace]
Let $m_\tau(p)=\tau(p)$, i.e., $f(t)=t$. Then $\delta_{m_\tau}=0$, and
\[
LL_{m_\tau}(R)=L_{m_\tau}(R)=S(R).
\]
The same conclusion holds for $m_\alpha(p)=\tau(p)^\alpha$ with
$0<\alpha\le1$.
\end{ex}

\begin{ex}[The indicator measure]
Let
\[
m_0(p)=\begin{cases}0,&p=0,\\1,&p\neq0.\end{cases}
\]
Here $f(t)=\e_{\{t>0\}}$ and $\delta_{m_0}=1$, so $LL_{m_0}(R)=L_{m_0}(R)=R\subsetneq S(R)$.
Likewise, $m(p)=\e_{\{p\neq0\}}+\tau(p)^\alpha$,
$0<\alpha\le1$, has $\delta_m=1$ and yields $LL_m(R)=R$.
\end{ex}

\begin{ex}[An infinite-valued measure]
Let $m_\infty(0)=0$ and $m_\infty(p)=\infty$ for $p\neq0$.
It satisfies Ciach's axioms, has $\delta_{m_\infty}=\infty$,
and $LL_{m_\infty}(R)=R$.
\end{ex}

\begin{rem} For a type $\mathrm{II}_1$ factor, the algebra $LL_m(R)$ depends on
$m$ only through whether $\delta_m$ vanishes. Thus no proper
intermediate algebra $R\subsetneq LL_m(R)\subsetneq S(R)$ arises in
the factorial case. Proper intermediate algebras may nevertheless
occur for type II von Neumann algebras with nontrivial center.
In the next stage, derivations may therefore be studied separately
on the two algebras $R$ and $S(R)$.
\end{rem}

%\begin{center}
%  {\Large\bfseries Stage III. Derivations on Algebras of Locally\par
%  $m$-Measurable Operators over Type $\mathrm{II}_1$ Factors}\par
%  \vspace{5pt}
%  {\small A rigorous classification and a constructive comparison of derivation spaces}
%\end{center}

\section{Derivations on algebras of locally $m$-measurable operators over type II$_1$ factors.}

Let $A$ be an algebra. A \emph{derivation} $D:A\to A$ is a complex-linear map satisfying
$D(xy)=D(x)y+xD(y)$. For any $a\in A$, set
\[
 ad_a(x)=[a,x]=ax-xa.
\]
The vector space $Der(A)$ is a complex Lie algebra under the commutator
$[D_1,D_2]=D_1D_2-D_2D_1$. We use the usual strong (closed-operator)
sum and product in $S(R)$. No continuity or $*$-preservation is assumed
for derivations.
%\section{The centers of $R$ and $S(R)$}
\begin{lem}[Relative commutant]\label{lem:relative_commutant}
If $t\in S(R)$ commutes with every element of $R$, then
$t\in \mathbb{C}\e$. Consequently,
\[
  Z(R)= \mathbb{C}\e,\qquad Z(S(R))= \mathbb{C}\e.
\]
In particular, if $a,b\in S(R)$ satisfy
$[a,r]=[b,r]$ for every $r\in R$, then $a-b\in \mathbb{C}\e$.
\end{lem}
\begin{proof}
Suppose that $[t,r]=0$ for all $r\in R$.
For each unitary $u\in R$ we have $utu^*=t$ in $S(R)$.
Taking adjoints and products gives
\[
  u(t^*t)u^*=t^*t,
  \qquad u|t|u^*=|t|.
\]
Therefore every spectral projection $E_B(|t|)$, for Borel sets
$B\subset[0,\infty)$, commutes with every unitary in $R$.
As $t$ is affiliated with $R$, these spectral projections also belong to
$R$. Thus
\[
  E_B(|t|)\in Z(R)= \mathbb{C}\e.
\]
Each such spectral projection is hence either $0$ or $1$.
By the spectral theorem and the dense definition of $|t|$, this implies
$|t|=\lambda1$ for some finite $\lambda\ge0$.
Consequently $t$ is bounded, and since $t$ is affiliated with $R$,
$t\in R$. As $t$ commutes with $R$ and $R$ is a factor,
$t\in Z(R)=\mathbb{C}\e$.

Every element of $Z(S(R))$ commutes with $R\subset S(R)$, and so lies in
$\mathbb{C}\e$ by the preceding argument. The reverse inclusion is immediate.
The final assertion follows by applying the first statement to $t=a-b$.
\end{proof}

We recall two known innerness results:
Sakai's theorem for derivations on von Neumann
algebras \cite{Sakai}, and the innerness theorem
for Murray--von Neumann algebras
\cite{BKS}.
\begin{thm}\label{thm:bounded}
Every complex-linear derivation $D:R\to R$ is inner and has an
implementing element in $R$. The implementer is unique modulo scalar
multiples of the identity. In particular,
\[
 Der(R)=Inn(R)\cong_{\mathrm{Lie}}R/ \mathbb{C}\e.
\]
\end{thm}
\begin{thm}\label{thm:unbounded}
Every complex-linear derivation $D:S(R)\to S(R)$ is inner, with an
implementing element belonging to $S(R)$. The implementer is unique
modulo scalar multiples of the identity. Thus
\[
Der(S(R))=Inn(S(R))\cong_{\mathrm{Lie}}S(R)/\mathbb{C}\e.
\]
No continuity assumption is necessary.
\end{thm}

\begin{thm}\label{thm:dichotomy_der}
Let $R$ be a type $\mathrm{II}_1$ factor, and let $m$ be any
subadditive measure on $ P(R)$. Then
\[
 Der(LL_m(R))=Inn(LL_m(R))
 \cong_{\mathrm{Lie}}
 \begin{cases}
 R/ \mathbb{C}\e,&\delta_m>0,\\[2pt]
 S(R)/ \mathbb{C}\e,&\delta_m=0.
 \end{cases}
\]
The first case includes $\delta_m=\infty$.
\end{thm}

\begin{proof}
By the above the algebra $LL_m(R)$ equals $R$ when $\delta_m>0$ and equals $S(R)$ when
$\delta_m=0$. Apply Theorems~\ref{thm:bounded} and
\ref{thm:unbounded}, respectively.
\end{proof}

\subsection*{A proper embedding of the two Lie algebras.}

\begin{prop}\label{prop:embedding}
There is a canonical injective Lie algebra homomorphism
\[
 \mathcal J: Der(R)\longrightarrow Der(S(R)),
 \qquad \mathcal J(ad_a|_R)=ad_a|_{S(R)}
 \quad(a\in R).
\]
It is not surjective. Consequently,
\[
\mathcal J(Der(R))\subsetneq Der(S(R)).
\]
Every derivation on $R$ admits exactly one extension to a derivation
on $S(R)$.
\end{prop}

\begin{proof}
By Theorem~\ref{thm:bounded}, each derivation on $R$ is $ad_a|_R$
for some $a\in R$. If $ad_a|_R=ad_b|_R$, then $a-b\in \mathbb{C}\e$,
and hence $ad_a|_{S(R)}=ad_b|_{S(R)}$.
Thus $\mathcal J$ is well-defined. It is evidently complex-linear and,
by the commutator identity for inner derivations, a Lie algebra
homomorphism. If $\mathcal J(ad_a|_R)=0$, then $a$ is central in
$S(R)$, hence scalar by Lemma~\ref{lem:relative_commutant}. Therefore
$ad_a|_R=0$, and $\mathcal J$ is injective.

We give an explicit derivation outside its image.
Using the diffuse trace scale of $R$, choose mutually orthogonal
projections $e_n\in R$ such that
\[
  \tau(e_n)=2^{-n},\qquad \bigvee_{n\ge1}e_n=\e.
\]
Let
\[
 h=\sum_{n=1}^{\infty}n e_n,
\]
where the sum is defined by the spectral theorem, and $h$ acts as
$n$ on $e_nH$. It is a densely defined, positive self-adjoint
operator affiliated with $R$. Hence $h\in S(R)$, but
$h\notin R$ since $\|h\|=\infty$.

The map $D_h=ad_h$ is a derivation on $S(R)$ by
Theorem~\ref{thm:unbounded}. Suppose that $D_h=\mathcal J(ad_a)$
for some $a\in R$. Then $[h-a,x]=0$ for every $x\in S(R)$, and in
particular for every $x\in R$.
By Lemma~\ref{lem:relative_commutant},
\[
 h-a\in \mathbb{C}\e.
\]
This would imply $h\in R$, a contradiction. Hence $D_h$ is not in the
image of $\mathcal J$, proving the inclusion is proper.

Finally, suppose $D\in Der(R)$ has another extension
$\widetilde D\in Der(S(R))$. By Theorem~\ref{thm:unbounded},
$\widetilde D= ad_k$ for some $k\in S(R)$, whereas
$D=ad_a|_R$ for some $a\in R$.
Since $[k-a,r]=0$ for every $r\in R$,
Lemma~\ref{lem:relative_commutant} implies $k-a\in \mathbb{C}\e$.
Thus $\widetilde D=ad_a|_{S(R)}=\mathcal J(D)$, giving uniqueness.
\end{proof}

\begin{prop}\label{prop:unitary-earlier}
For the explicit operator $h=\sum_{n\ge1}ne_n$ constructed above,
there exists a self-adjoint unitary $u\in R$ such that $[h,u]\in S(R)\setminus R$.
In particular, $ad_h$ does not map $R$ into $R$.
\end{prop}

\begin{proof}
Since $\tau(e_1)=1/2$ and $e_1Re_1$ is diffuse, choose mutually
orthogonal projections $q_n\le e_1$ satisfying
\[
 \tau(q_n)=2^{-(n+1)},\qquad
 \bigvee_{n\ge1}q_n=e_1.
\]
For every $n\ge1$ we have
$\tau(q_n)=\tau(e_{n+1})$, and therefore there exists a partial
isometry $v_n\in R$ with $v_n^*v_n=q_n$, \ $v_nv_n^*=e_{n+1}$.
The supports of the $v_n$ are mutually orthogonal. The strong sum
$v=\sum_{n\ge1}v_n$ exists in $R$ and satisfies
$v^*v=e_1$, \ $vv^*=\e-e_1$, \ $v^2=0$.
Hence $u=v+v^*$ is a self-adjoint unitary in $R$.
For $n\ge1$ and a unit vector $\xi\in q_nH$, we have
\[
 h\xi=\xi,\qquad u\xi=v_n\xi\in e_{n+1}H,
 \qquad h(u\xi)=(n+1)u\xi.
\]
It follows that $[h,u]\xi=hu\xi-uh\xi=n\,u\xi$, \ $\|[h,u]\xi\|=n$.
Since $n$ is arbitrary, $[h,u]$ is unbounded.
The commutator belongs to $S(R)$ because $S(R)$ is a $*$-algebra
containing both $h$ and $u$. Hence
$[h,u]\in S(R)\setminus R$.
\end{proof}

\begin{cor}\label{cor:relative-earlier}
Let $Der(R,S(R))$ denote the complex vector space of complex-linear
derivations from $R$ into the natural $R$-bimodule $S(R)$.
Then
\[
 Der(R,S(R))\cong_{C}S(R)/ \mathbb{C}\e.
\]
Under this identification, the subspace $Der(R,R)$ corresponds to
$R/ \mathbb{C}\e$. Consequently there is a canonical vector-space isomorphism
\[
Der(R,S(R))/ Der(R,R)\cong_{C}S(R)/R.
\]
The quotient here is a quotient of complex vector spaces; no Lie
algebra structure on $Der(R,S(R))$ is asserted.
\end{cor}

\begin{proof}
Theorem~1.1 of Huang--Kudaybergenov--Sukochev~\cite{HKS}, app\-lied
to the $*$-subalgebra $R\subseteq S(R)$, states that every derivation
$D:R\to S(R)$ is spatially implemented by some $h\in S(R)$.
The map $h\mapsto ad_h|_R$ is therefore surjective and has kernel
$\mathbb{C}\e$ by Lemma~\ref{lem:relative_commutant}.
This proves the first vector-space isomorphism.
By Sakai's theorem, the subspace of derivations with range in $R$ is
precisely the set implemented by elements of $R$. Since
$\mathbb{C}\e\subseteq R\subseteq S(R)$, the quotient-of-quotients
isomorphism yields the final assertion.
\end{proof}
%\section{Conclusion and scope}
\begin{rem}
The innerness of derivations on $R$ and $S(R)$ is not asserted here
as a new result: it follows from the cited theorems of Sakai and of
Ber--Kudaybergenov--Sukochev. The new notation $LL_m(R)$ merely
selects one of these two known algebras in the type $\mathrm{II}_1$
factor setting. The explicit unitary construction displays a concrete
difference between the corresponding derivation spaces. It should not
be presented as a previously unknown theorem without a separate
novelty search.
\end{rem}

\section{Comparison of derivations associated with \\
different subadditive measures.}

We use Sakai's innerness theorem for derivations $R\to R$~\cite{Sakai}
and the Ber--Kudaybergenov--Sukochev innerness theorem for derivations
$S(R)\to S(R)$~\cite{BKS}. In particular,
\begin{equation}\label{eq:inner}
 Der(R)=\{ad_a|_R:a\in R\},\qquad
 Der(S(R))=\{ad_b|_{S(R)}:b\in S(R)\}.
\end{equation}
No continuity hypothesis is imposed on these complex-linear derivations.
%\section{Comparing two subadditive measures}
\begin{thm}[Complete comparison criterion]\label{thm:comparison}
Let $m_1,m_2$ be subadditive measures on $P(R)$, and put
\[
 A_i=LL_{m_i}(R),\qquad
 \delta_i=\inf\{m_i(p):p\in P(R),\ p\ne0\},
 \quad i=1,2.
\]
Then:
\begin{enumerate}
\item $A_1=A_2$ if and only if either $\delta_1,\delta_2>0$,
      or $\delta_1=\delta_2=0$.
\item $A_1\subsetneq A_2$ if and only if $\delta_1>0$ and
      $\delta_2=0$.
\item $A_1\subseteq A_2$ if and only if
      $\delta_1>0$ or $\delta_2=0$.
\end{enumerate}
In particular, the two algebras $A_1,A_2$ are always comparable.
\end{thm}
\begin{proof}
By Theorem \ref{thm:dichotomy}, each $A_i$ equals either $R$ or $S(R)$,
with the alternatives determined exactly by $\delta_i$.
To see that the two algebras are distinct, choose pairwise orthogonal
nonzero projections $e_n\in R$ with $\tau(e_n)=2^{-n}$.
The operator
\[
 h=\sum_{n=1}^\infty ne_n
\]
is densely defined, self-adjoint and affiliated with $R$, so it lies
in $S(R)$; it is unbounded, hence does not belong to $R$.
Thus $R\subsetneq S(R)$. All three assertions now follow by examining
the four possible pairs $(R,R)$, $(R,S(R))$, $(S(R),R)$, and
$(S(R),S(R))$.
\end{proof}

\begin{rem}
The criterion depends only on the dichotomy parameters $\delta_i$,
not on a pointwise comparison between $m_1$ and $m_2$. In particular,
there is no general need to assume $m_1(p)\le m_2(p)$.
\end{rem}

\begin{thm}[Which derivations of $S(R)$ preserve $R$?]
\label{thm:stabilizer}
For a derivation $\Delta\in Der(S(R))$, the following conditions
are equivalent:
\begin{enumerate}
\item $\Delta(R)\subseteq R$;
\item $\Delta=\mathcal J(D)$ for some (necessarily unique)
$D\in Der(R)$;
\item $\Delta= ad_a|_{S(R)}$ for some $a\in R$.
\end{enumerate}
Equivalently,
\[
\mathcal J(Der(R))=
 \{\Delta\in Der(S(R)):\Delta(R)\subseteq R\}.
\]
\end{thm}
\begin{proof}
Suppose $\Delta(R)\subseteq R$. Then its restriction
$D=\Delta|_R$ is a derivation $R\to R$. By Proposition~\ref{prop:embedding},
$\mathcal J(D)$ is the unique extension of $D$ to $S(R)$. Since
$\Delta$ is another such extension, $\Delta=\mathcal J(D)$.
This proves (1)$\Rightarrow$(2).

Condition (2) implies (3) by the definition of $\mathcal J$.
Finally, if $\Delta=ad_a$ with $a\in R$, then
$[a,r]\in R$ for every $r\in R$, proving (3)$\Rightarrow$(1).
\end{proof}

\begin{cor}[Comparison of derivation algebras]\label{cor:measure_maps}
Suppose $A_1\subseteq A_2$ in Theorem~\ref{thm:comparison}.
If $A_1=A_2$, their derivation Lie algebras coincide canonically.
If $A_1\subsetneq A_2$, then $A_1=R$, $A_2=S(R)$, and the unique
canonical Lie algebra embedding
\[
 Der(A_1)\longrightarrow Der(A_2)
\]
obtained by extending each derivation is $\mathcal J$.
Its image consists precisely of those derivations of $A_2$
which preserve $A_1$.
\end{cor}
\begin{proof}
Apply Theorems~\ref{thm:comparison}, \ref{thm:stabilizer}, and Proposition \ref{prop:embedding}.
\end{proof}

\subsection*{A proper non-ideal Lie subalgebra.}

\begin{thm}[Properness and failure of ideality]\label{thm:nonideal}
The canonical Lie subalgebra $\mathcal J(Der(R))$ of
$Der(S(R))$ is proper and is not a Lie ideal.
In particular, there is no induced Lie algebra structure on the
ordinary vector-space quotient
$Der(S(R))/\mathcal J(Der(R))$ coming from a quotient
Lie algebra construction.
\end{thm}
\begin{proof}
Let $h,u$ be as in Proposition~\ref{prop:unitary-earlier}.
If $ad_h\in\mathcal J(Der(R))$, there would exist $a\in R$
with $ad_h=ad_a$ on $S(R)$. By Lemma~\ref{lem:relative_commutant},
$h-a\in \mathbb{C}\e$, contradicting $h\notin R$. Thus the image is proper.

Since $u\in R$, $ad_u\in\mathcal J(Der(R))$. However,
\[
 [ad_h, ad_u]=ad_{[h,u]}.
\]
If $ad_{[h,u]}$ were in $\mathcal J(Der(R))$, there would be
$a\in R$ with $[h,u]-a\in \mathbb{C}\e$ by Lemma~\ref{lem:relative_commutant}.
This would force $[h,u]\in R$, contrary to
Proposition~\ref{prop:unitary-earlier}. Hence the image is not a Lie ideal.
\end{proof}
%\section{The quotient of derivation spaces}

\begin{thm}[Vector-space quotient and infinite dimension]
\label{thm:quotient}
There is a canonical isomorphism of complex vector spaces
\[
 Der(S(R))/\mathcal J(Der(R))
 \ \cong_{C}\ S(R)/R.
\]
Both quotient vector spaces are infinite-dimensional.
They are not asserted to be quotients of Lie algebras.
\end{thm}
\begin{proof}
By~\eqref{eq:inner} and Lemma~\ref{lem:relative_commutant},
$\Phi:S(R)\longrightarrow Der(S(R))$, \ $\Phi(b)= ad_b$
is a surjective complex-linear map with kernel $\mathbb{C}\e$.
Moreover, $\Phi(R)=\mathcal J(Der(R))$.
Since $\mathbb{C}\e\subseteq R$, $\Phi$ induces a well-defined surjective
linear map
\[
 \overline\Phi:S(R)/R\longrightarrow
 Der(S(R))/\mathcal J(Der(R)),\quad
 b+R\longmapsto ad_b+\mathcal J(Der(R)).
\]
If $\overline\Phi(b+R)=0$, then $ad_b=ad_a$ on $S(R)$
for some $a\in R$, so $b-a\in \mathbb{C}\e$ by Lemma~\ref{lem:relative_commutant}.
It follows that $b\in R$. Hence $\overline\Phi$ is injective.

To prove infinite dimension, retain the orthogonal projections
$e_n$ from Proposition~\ref{prop:unitary-earlier}, and for each integer
$k\ge1$ set
\[
 h_k=\sum_{n=1}^\infty n^k e_n\in S(R).
\]
Suppose a finite linear combination
$\displaystyle\sum_{k=1}^{N}\alpha_kh_k$ belongs to $R$. On $e_nH$ it acts as
$P(n)1$, where $\displaystyle P(t)=\sum_{k=1}^{N}\alpha_kt^k$.
Boundedness of this operator implies that
$\displaystyle\sup_{n\ge1}|P(n)|<\infty$. A nonzero polynomial with zero constant
term cannot be bounded along the positive integers. Thus
$\alpha_1=\cdots=\alpha_N=0$.
Consequently the cosets $h_k+R$ are linearly independent in $S(R)/R$.
Therefore both quotient vector spaces are infinite-dimensional.
\end{proof}

\begin{rem}[An explicit pair of measures]
For
\[
 m_0(p)=\e_{\{p\ne0\}},\qquad m_\tau(p)=\tau(p),
\]
both of which satisfy Ciach's axioms, one has
\[
 \delta_{m_0}=1,\quad \delta_{m_\tau}=0,\qquad
 LL_{m_0}(R)=R\subsetneq S(R)=LL_{m_\tau}(R).
\]
The resulting comparison of derivation spaces
is described by Proposition~\ref{prop:embedding}
and Theorems~\ref{thm:stabilizer},
\ref{thm:nonideal}, and \ref{thm:quotient}.
For $m_\alpha(p)=\tau(p)^\alpha$, $0<\alpha\le1$, the value of
$\delta_{m_\alpha}$ is also zero and the algebra remains $S(R)$.
\end{rem}

\begin{rem}[Interpretation and novelty]
The classification of $LL_m(R)$ into $R$ or $S(R)$ is the outcome of
the previously established type $\mathrm{II}_1$ dichotomy.
The innerness theorems for $R$ and $S(R)$ are due to Sakai and to
Ber--Kuday\-bergenov--Sukochev. The comparisons developed here are
formal consequences of those results, together with explicit
constructions. They should not be advertised as previously unknown
without a separate exhaustive literature review. To seek new results,
one should next investigate nonfactorial type $\mathrm{II}$ algebras or
value-restricted derivations and other additional structures.
\end{rem}

\section{Range-restricted derivations on locally $m$-measurable operator algebras.}

Define $A_m=LL_m(R)$. By the above, we have
\begin{equation}\label{eq:dichotomy}
 A_m=L_m(R)=
 \begin{cases}
 R,&\delta_m>0,\\
 S(R),&\delta_m=0.
 \end{cases}
\end{equation}
The first alternative includes $\delta_m=+\infty$.
The equality $LL_m(R)=L_m(R)$ follows directly from
$P(Z(R))=\{0,1\}$.

For a unital subalgebra $A\subseteq S(R)$ and any complex linear
subspace $B\subseteq S(R)$, define the \emph{range-
restricted derivation space}
\begin{equation}\label{eq:relative-definition}
\hspace*{-0.7cm}
Der_{S(R)}(A;B)=\left\{D:A\to S(R):
\begin{array}{l}
 D\text{ is complex-linear},\\
 D(xy)=D(x)y+xD(y), x,y\in A\\
 D(A)\subseteq B
 \end{array}\right\}
\end{equation}
The Leibniz identity is an identity in $S(R)$.
\emph{We do not require $B$ to be an $A$-bimodule}; consequently
\eqref{eq:relative-definition} makes sense even for $A=S(R)$ and
$B=R$. These sets are complex vector spaces, but need not carry a
natural Lie bracket under composition if their values do not lie
in their domains. We abbreviate $Der(A;B)=Der_{S(R)}(A;B)$
in the sequel.

We shall use the following established results:
\begin{enumerate}
\item[(i)] (Sakai) Every derivation $R\to R$ is inner, implemented
by an element of $R$~\cite{Sakai}.
\item[(ii)] (Ber--Kudaybergenov--Sukochev) Every complex-linear
derivation $S(R)\to S(R)$ is inner~\cite{BKS}.
\item[(iii)] (Huang--Kudaybergenov--Sukochev) Every complex-linear
derivation $R\to S(R)$ is of the form $ad_a|_R$ for some
$a\in S(R)$~\cite[Theorem~1.1]{HKS}.
\end{enumerate}
No continuity assumptions are imposed on any of our derivations.

The next lemma is proved from the structure of a diffuse finite
factor. In particular, it does not rely on the previously published
proper-bimodule rigidity theorem.

\begin{lem}[Unbounded commutator witness]\label{lem:witness}
For every $a\in S(R)\setminus \mathbb{C}\e$, there exists a positive,
unbounded $x\in S(R)$ such that
\[
 [a,x]\in S(R)\setminus R.
\]
In fact, $x$ may be chosen supported on a nonzero projection
$e\in R$ and with a purely discrete spectral decomposition
$\displaystyle x=\sum_{n\ge1}n e_n$, where $e_n\in P(R)$ are nonzero and
pairwise orthogonal.
\end{lem}
\begin{proof}
Since $a$ does not commute with all of $R$ by
Lemma~\ref{lem:relative_commutant}, it does not commute with all the
projections of $R$. We first claim that there exists
$p\in P(R)$ with
\begin{equation}\label{eq:offdiagonal}
 b:=pa(\e-p)\ne0.
\end{equation}
Indeed, if $pa(\e-p)=0$ for every projection $p$, applying the same
assumption to $\e-p$ gives $(\e-p)ap=0$. It follows that
$pa=pap=ap$ for every projection $p$. By the spectral theorem,
linear combinations of projections approximate any self-adjoint
element of $R$ in operator norm, hence in the measure topology.
Since left and right multiplication by a fixed element of $S(R)$
are continuous in the measure topology, $a$ commutes with all
of $R$, a contradiction. Put $q=\e-p$. Then $b=paq$ is a
nonzero affiliated operator satisfying $b=pbq$.

Consider the spectral projections of the positive affiliated
operator $|b|=(b^*b)^{1/2}$. Because $b\ne0$, there exist finite
constants $0<\varepsilon<K<\infty$ for which
\begin{equation}\label{eq:e}
 e=E_{[\varepsilon,K]}(|b|)\ne0.
\end{equation}
As $b=bq$, we have $e\le q$. Spectral calculus yields
\begin{equation}\label{eq:lowerbound}
 \varepsilon\|\xi\|\le\|b\xi\|\le K\|\xi\|
 \qquad(\xi\in eH).
\end{equation}
In particular $be$ is a bounded operator affiliated with $R$,
so $be\in R$.

The nonzero corner $eRe$ is a finite diffuse von Neumann algebra
(in fact, a type $\mathrm{II}_1$ factor with normalized trace
$\tau(\cdot)/\tau(e)$). By repeated trace splitting, choose
pairwise orthogonal nonzero projections $e_n\le e$ such that
\[
 \tau(e_n)=2^{-n}\tau(e),\qquad n\ge1,
 \qquad \bigvee_{n\ge1}e_n=e.
\]
Define the positive operator
\begin{equation}\label{eq:x}
 x=\sum_{n=1}^{\infty}n e_n.
\end{equation}
More precisely, this is the spectral sum on the domain
\[
 \mathcal D(x)=(\e-e)H\oplus
 \left\{\xi\in eH:\sum_{n\ge1}n^2\|e_n\xi\|^2<\infty\right\}.
\]
This domain contains $(\e-e)H$ and every finite direct sum
$\displaystyle\bigoplus_{n=1}^N e_nH$, so it is dense in $H$.
The spectral theorem shows that $x$ is positive self-adjoint and
affiliated with $R$, hence $x\in S(R)$. It is unbounded because
$xe_n=ne_n$ and all $e_n$ are nonzero. We also have
$x=exe=qxq$ and $px=0$, with the products interpreted in $S(R)$.

Using the associative algebra operations in $S(R)$,
for every $n\ge1$ we have
\begin{align}\label{eq:corner}
 p[a,x]e_n
 &=p(ax-xa)e_n\\
 &=pa(xe_n)-(px)ae_n\notag\\
 &=n\,pae_n=n\,be_n.\notag
\end{align}
Notice that $be_n=(be)e_n\in R$, and inequality
\eqref{eq:lowerbound} gives
\[
 \|be_n\|\ge\varepsilon
 \qquad(n\ge1),
\]
since each $e_n\ne0$. Hence
\begin{equation}\label{eq:normgrowth}
 \|p[a,x]e_n\|\ge n\varepsilon.
\end{equation}
If the closed affiliated commutator $[a,x]$ belonged to $R$, then
$\|p[a,x]e_n\|\le\|[a,x]\|$ for every $n$, contradicting
\eqref{eq:normgrowth}. Therefore $[a,x]\notin R$.
\end{proof}

\begin{thm}[Vanishing of derivations with values in $R$]
\label{thm:vanishing}
Let $R$ be any type $\mathrm{II}_1$ factor. If
$D:S(R)\to S(R)$ is a complex-linear derivation such that
\[
 D(S(R))\subseteq R,
\]
then $D=0$. Equivalently,
\begin{equation}\label{eq:zero}
Der(S(R);R)=\{0\}.
\end{equation}
\end{thm}
\begin{proof}
By the innerness theorem for Murray--von Neumann algebras
\cite{BKS}, there exists $a\in S(R)$ with
$D=ad_a$. If $a\notin \mathbb{C}\e$, Lemma~\ref{lem:witness} produces an
$x\in S(R)$ such that $D(x)=[a,x]\notin R$, contradicting the range
assumption. Thus $a\in \mathbb{C}\e$, whence $D=0$.
\end{proof}

\begin{rem}[Source and historical priority]\label{rem:priority}
The conclusion of Theorem~\ref{thm:vanishing} is \emph{not new}.
The special case with values in $R$ was already established
by Kadison and Liu~\cite{KadisonLiu} in 2014.
It also follows immediately from the stronger theorem of
Ber--Kudaybergenov--Sukochev~\cite{BKS2}, which states that
any derivation from $S(R)$ into a proper $R$-bimodule of $S(R)$
is zero. The proof above is included because it gives an explicit
positive unbounded witness for every noncentral implementer and
isolates the mechanism needed for the two-measure classification.
\end{rem}
%\section{All range restrictions associated with two measures}
\begin{thm}[Complete two-measure range classification]
\label{thm:two-measure}
Let $m,n$ be arbitrary Ciach subadditive measures on $P(R)$.
Put
\[
 A_m=LL_m(R),\qquad B_n=L_n(R),
 \qquad \delta_m=\inf_{p\ne0}m(p),\quad
 \delta_n=\inf_{p\ne0}n(p).
\]
Then, inside $S(R)$, the following exhaustive classification holds:
\begin{equation}\label{eq:fourcases}
Der(A_m;B_n)=
\begin{cases}
 \{ad_a|_R:a\in R\},
     &\delta_m>0,\ \delta_n>0,\\[2pt]
 \{ad_a|_R:a\in S(R)\},
     &\delta_m>0,\ \delta_n=0,\\[2pt]
 \{0\},
     &\delta_m=0,\ \delta_n>0,\\[2pt]
 \{ad_a|_{S(R)}:a\in S(R)\},
     &\delta_m=0,\ \delta_n=0.
\end{cases}
\end{equation}
In each nonzero case, the implementing element is unique modulo
$\mathbb{C}\e$. All four cases include measures taking the value $+\infty$.
\end{thm}
\begin{proof}
By Theorem \ref{thm:dichotomy} $A_m\in\{R,S(R)\}$, \ $B_n\in\{R,S(R)\}$,
and the alternatives are determined by $\delta_m$ and
$\delta_n$. Thus there are exactly four possibilities.

If $A_m=B_n=R$, the assertion follows from Sakai's innerness
theorem~\cite{Sakai}. If $A_m=R$ and $B_n=S(R)$, it follows from
the theorem of Huang--Kudayber\-genov--Sukochev
\cite[Theorem~1.1]{HKS}. If $A_m=S(R)$ and $B_n=R$, it follows
from Theorem~\ref{thm:vanishing}. Finally, if
$A_m=B_n=S(R)$, it follows from the innerness theorem
\cite{BKS}. In each nonzero case, uniqueness of implementers
modulo $\mathbb{C}\e$ is Lemma~\ref{lem:relative_commutant}.
\end{proof}

\begin{cor}[Range in $R$ for one measure]\label{cor:one-measure}
For each subadditive measure $m$, one has
\[
 Der(LL_m(R);R)=
 \begin{cases}
 \{ad_a|_R:a\in R\},&\delta_m>0,\\
 \{0\},&\delta_m=0.
\end{cases}
\]
On the other hand, the restriction $D(LL_m(R))\subseteq L_m(R)$
is automatically satisfied by every derivation
$D:LL_m(R)\to LL_m(R)$, because $LL_m(R)=L_m(R)$ for a factor.
Consequently,
\[
Der(LL_m(R);L_m(R))=Der(LL_m(R)).
\]
\end{cor}
\begin{proof}
The first statement is the first and third rows of
Theorem~\ref{thm:two-measure}, since $R$ is the target.
The second is the equality $LL_m(R)=L_m(R)$.
\end{proof}

\begin{cor}[Canonical restriction and extension]
\label{cor:restriction}
Restriction defines an isomorphism of complex vector spaces
\[
 \operatorname{res}:Der(S(R);S(R))
 \stackrel{\cong}{\longrightarrow} Der(R;S(R)),
 \quad \Delta\longmapsto\Delta|_R.
\]
In particular, every complex-linear derivation $R\to S(R)$
has a unique extension to a derivation $S(R)\to S(R)$.
Moreover, for $a\in S(R)$,
\[
 ad_a(R)\subseteq R
 \quad\Longleftrightarrow\quad
 a\in R+\mathbb{C}\e=R.
\]
\end{cor}
\begin{proof}
Every derivation $D:R\to S(R)$ equals $ad_a|_R$ for some
$a\in S(R)$ by~\cite[Theorem~1.1]{HKS}, and $ad_a$ defines
a derivation of $S(R)$. Hence restriction is surjective.
If $ad_a|_R=0$ for some $a\in S(R)$, then
$a\in \mathbb{C}\e$ by Lemma~\ref{lem:relative_commutant}; therefore $ad_a=0$
on all of $S(R)$. By innerness of all derivations on $S(R)$,
restriction is also injective.

If $ad_a(R)\subseteq R$, then its restriction to $R$ is a
derivation $R\to R$. By Sakai's theorem, it equals
$ad_b|_R$ for some $b\in R$. The relative-commutant lemma gives
$a-b\in \mathbb{C}\e$, so $a\in R+\mathbb{C}\e=R$. The converse follows from
$[R,R]\subseteq R$.
\end{proof}
\begin{rem}
It is easy to see that Corollary \ref{cor:restriction} strengthens the statement of Proposition \ref{prop:embedding}.
\end{rem}
\begin{thm}[Ber--Kudaybergenov--Sukochev, 2020]
\label{thm:published}
Let $R$ be a type $\mathrm{II}_1$ factor and let
$\mathcal B\subsetneq S(R)$ be a complex linear $R$-bimodule;
that is,
\[
 R\mathcal B\subseteq\mathcal B,
 \qquad \mathcal B R\subseteq\mathcal B.
\]
Then any derivation $D:S(R)\to S(R)$ satisfying
$D(S(R))\subseteq\mathcal B$ is identically zero.
\end{thm}
\begin{cor}[Classical noncommutative $L^p$ targets]
\label{cor:lp}
Let $0<p\le\infty$, and let $L^p(R,\tau)$ denote the associated
noncommutative $L^p$ space (where $L^\infty(R,\tau)=R$).
Then
\[
Der(S(R);L^p(R,\tau))=\{0\}.
\]
\end{cor}
\begin{proof}
For every $0<p\le\infty$, the space $L^p(R,\tau)$ is a
complex linear $R$-bimodule of $S(R)$. It is proper:
for $p=\infty$, this follows from $R\subsetneq S(R)$.
For $0<p<\infty$, choose pairwise orthogonal nonzero projections
$e_n$ with $\tau(e_n)=2^{-n}$ and let
\[
 y=\sum_{n\ge1}2^{n/p}e_n\in S(R).
\]
Then
\[
 \tau(|y|^p)=\sum_{n\ge1}2^n\tau(e_n)=\infty,
\]
so $y\notin L^p(R,\tau)$. The conclusion therefore follows from
Theorem~\ref{thm:published}.
\end{proof}
%\section{Explicit measures and mathematical interpretation}
\begin{prop}[A two-measure example]\label{prop:example}
Let
$m_0(p)=\e_{\{p\ne0\}}$, \ $m_\tau(p)=\tau(p)$.
Then
$LL_{m_0}(R)=L_{m_0}(R)=R$, \ $LL_{m_\tau}(R)=L_{m_\tau}(R)=S(R)$,
and the four relative derivation spaces are
\begin{align*}
 Der(R;R)&=\{ad_a|_R:a\in R\},\\
 Der(R;S(R))&=\{ad_a|_R:a\in S(R)\},\\
 Der(S(R);R)&=\{0\},\\
 Der(S(R);S(R))&=\{ad_a|_{S(R)}:a\in S(R)\}.
\end{align*}
The second space is strictly larger than the first (under their
natural identi\-fi\-cation as derivations from $R$ into $S(R)$).
\end{prop}
\begin{proof}
The equalities of operator algebras follow from
$\delta_{m_0}=1$ and $\delta_{m_\tau}=0$. The four formulae are
Theorem~\ref{thm:two-measure}.
For strictness choose $a\in S(R)\setminus R$, which exists because
$R$ is diffuse. If $ad_a|_R=ad_b|_R$ for some $b\in R$, then
$a-b\in \mathbb{C}\e$ by Lemma~\ref{lem:relative_commutant}, contradicting
$a\notin R$.
\end{proof}

\begin{rem}[Scope and novelty]\label{rem:scope}
The vanishing of derivations $S(R)\to R$ and, more generally,
$S(R)\to\mathcal B$ for every proper $R$-bimodule
$\mathcal B\subsetneq S(R)$ had already been established in
\cite{BKS2}. Thus Theorem~\ref{thm:vanishing} is an alternative
constructive proof of a known special case, not an original
innerness or vanishing result. The four-case table in
Theorem~\ref{thm:two-measure} follows by combining that established
rigidity with~\eqref{eq:dichotomy} and known innerness theorems.
For the search for genuinely new results, one should examine
nonfactorial type II von Neumann algebras, where the algebras
$LL_m(M)$ need not reduce to $M$ or $S(M)$, and where proper
$M$-bimodule restrictions can interact nontrivially with central
localization.
\end{rem}

\section{Measure topologies and automatic continuity for locally $m$-measurable operator algebras.}

For $x\in S(R)$ and $\varepsilon>0$ put
\begin{equation}\label{eq:tail}
 q_\varepsilon(x):=E_{(\varepsilon,\infty)}(|x|)\in P(R).
\end{equation}
Following Ciach's Section~2, for $\varepsilon,\eta>0$ define
\begin{equation}\label{eq:ciach-neighborhood}
 V_m(\varepsilon,\eta):=\left\{x\in L_m(R):
 \begin{array}{l}
 \exists\,p\in P(R)\text{ with }pH\subseteq\mathcal D(x),\\[-1pt]
 \|xp\|\le\varepsilon,\quad m(1-p)\le\eta
 \end{array}\right\}.
\end{equation}
These sets generate the $m$-measure topology $t_m$ on $L_m(R)$.
As Ciach proves, $t_m$ makes $L_m(R)$ a complete metrizable
topological $*$-algebra
\cite[Corollaries~2.5--2.7]{Ciach}. In particular the product
$L_m(R)\times L_m(R)\to L_m(R)$ is jointly continuous.
For $m=\tau$, denote this topology by $t_\tau$, the usual
measure topology on $S(R)$ \cite{Nelson}.
%\section{Spectral neighborhoods and the topological dichotomy}
\begin{lem}[Exact spectral description of neighborhoods]
\label{lem:spectral}
For every Ciach subadditive measure $m$ on a type $\mathrm{II}_1$
factor, every $\varepsilon,\eta>0$, and every $x\in L_m(R)$,
\begin{equation}\label{eq:spectral-neighborhood}
 x\in V_m(\varepsilon,\eta)
 \quad\Longleftrightarrow\quad
 m\bigl(q_\varepsilon(x)\bigr)\le\eta.
\end{equation}
Consequently, a sequence $x_n\in L_m(R)$ converges to zero in
$t_m$ if and only if
\begin{equation}\label{eq:convergence}
 m\bigl(q_\varepsilon(x_n)\bigr)\longrightarrow0
 \quad\text{for every }\varepsilon>0.
\end{equation}
\end{lem}
\begin{proof}
Suppose $x\in V_m(\varepsilon,\eta)$, and choose $p$ as in
\eqref{eq:ciach-neighborhood}. Set $q=q_\varepsilon(x)$.
We claim that
\begin{equation}\label{eq:trace-dom}
 \tau(q)\le\tau(\e-p).
\end{equation}
Otherwise, using the dimension formula for projections in the
finite factor, we obtain
$\tau(q\wedge p)\ge\tau(q)+\tau(p)-1>0$.
In particular, $q\wedge p\ne0$; note that no commutativity
assumption on $q$ and $p$ is needed for the dimension formula.
Choose $0\neq\xi\in(q\wedge p)H$. Then $\xi\in pH\subseteq
\mathcal D(x)$ and, by the strict spectral cutoff defining $q$,
$\|x\xi\|=\||x|\xi\|>\varepsilon\|\xi\|$.
This contradicts $\|xp\|\le\varepsilon$. Thus
\eqref{eq:trace-dom} holds. By comparison of projections in
$R$, there is $r\le\e-p$ with $q\sim_R r$. Hence
$m(q)=m(r)\le m(\e-p)\le\eta$.
Conversely, if $m(q_\varepsilon(x))\le\eta$, take
$p=E_{[0,\varepsilon]}(|x|)=1-q_\varepsilon(x)$.
Then $pH\subseteq\mathcal D(x)$, $\|xp\|\le\varepsilon$,
and $m(\e-p)\le\eta$. Thus $x\in V_m(\varepsilon,\eta)$.
This proves \eqref{eq:spectral-neighborhood}.
Finally, \eqref{eq:convergence} follows directly from the
neighborhood basis: for every fixed $\varepsilon>0$, every
$\eta>0$ must eventually dominate $m(q_\varepsilon(x_n))$.
\end{proof}

\begin{thm}[Complete classification of the $m$-measure topology]
\label{thm:topology}
Let $R$ be a type II$_1$ factor, and let
$m$ be a subadditive measure. With $\delta_m$ as above,
\begin{equation}\label{eq:main-topology}
 (LL_m(R),t_m)=
 \begin{cases}
 (R,t_{\|\cdot\|}),&\delta_m>0,\\[2pt]
 (S(R),t_\tau),&\delta_m=0.
 \end{cases}
\end{equation}
The identifications are equality of the underlying operator algebras
and equality of their topologies, not merely abstract topological
isomorphisms. They hold without any separability or finite-valuedness
assumption on $m$.
\end{thm}
\begin{proof}
The underlying algebra equalities were established in
\eqref{eq:dichotomy}.

\smallskip\noindent
\emph{Case 1: $\delta_m>0$.}
Choose any finite number $0<\eta_0<\delta_m$. Such a number exists
also when $\delta_m=\infty$. Every nonzero projection $q$ satisfies
$m(q)\ge\delta_m>\eta_0$; therefore, for any $x\in R$,
\[
 m(q_\varepsilon(x))\le\eta_0
 \quad\Longleftrightarrow\quad q_\varepsilon(x)=0
 \quad\Longleftrightarrow\quad\|x\|\le\varepsilon.
\]
By Lemma~\ref{lem:spectral},
\begin{equation}\label{eq:norm-balls}
 V_m(\varepsilon,\eta_0)=\{x\in R:\|x\|\le\varepsilon\}.
\end{equation}
Moreover, for any $\eta>0$,
$\{\|x\|\le\varepsilon\}\subseteq V_m(\varepsilon,\eta)$.
Thus each $m$-neighborhood contains a norm neighborhood, while
\eqref{eq:norm-balls} provides arbitrarily small norm neighborhoods
among the $m$-neighborhoods. The two topologies coincide.

\smallskip\noindent
\emph{Case 2: $\delta_m=0$.}
Here $L_m(R)=S(R)$ and, by $m(p)=f(\tau(p))$, $f(0)=0$, $f(t)>0$  ($t>0$), and
 $\delta_m=\lim_{t\downarrow0}f(t)$, we obtain
\begin{equation}\label{eq:limitf}
 \lim_{s\downarrow0}f(s)=0,\qquad f(s)>0\quad(s>0).
\end{equation}
Fix $\varepsilon,\eta>0$. Choose $s\in(0,1]$ with
$f(s)<\eta$. If $x\in V_\tau(\varepsilon,s)$, then the spectral
lemma applied to $\tau$ gives
$\tau(q_\varepsilon(x))\le s$. Consequently,
\[
 m(q_\varepsilon(x))
 =f(\tau(q_\varepsilon(x)))\le f(s)<\eta.
\]
Thus
\begin{equation}\label{eq:inclusion-one}
 V_\tau(\varepsilon,s)\subseteq V_m(\varepsilon,\eta).
\end{equation}
Conversely, fix $\varepsilon>0$ and $s\in(0,1]$, and choose
$0<\eta<f(s)$, which is possible by \eqref{eq:limitf}.
If $x\in V_m(\varepsilon,\eta)$ and
$\tau(q_\varepsilon(x))\ge s$, monotonicity of $f$ gives
\[
 m(q_\varepsilon(x))
 =f(\tau(q_\varepsilon(x)))\ge f(s)>\eta,
\]
a contradiction. Hence $\tau(q_\varepsilon(x))<s$, and
\begin{equation}\label{eq:inclusion-two}
 V_m(\varepsilon,\eta)\subseteq V_\tau(\varepsilon,s).
\end{equation}
For a $\tau$-neighborhood with parameter $s>1$, one may replace $s$
by a number in $(0,1]$. Thus \eqref{eq:inclusion-one} and
\eqref{eq:inclusion-two} establish mutual cofinality of the two
neighborhood bases and prove $t_m=t_\tau$.
\end{proof}

\begin{cor}[Independence of the particular measure]
\label{cor:independence}
If $m_1,m_2$ are subadditive measures on $R$ with
$\delta_{m_1}$ and $\delta_{m_2}$ both positive or both zero,
then $LL_{m_1}(R)=LL_{m_2}(R)$ and the corresponding
measure topologies are equal. Consequently the pair
$(LL_m(R),t_m)$ depends only on whether $\delta_m$ vanishes.
\end{cor}
\begin{proof}
Apply Theorem~\ref{thm:topology} twice.
\end{proof}

\begin{rem}[The associated metric]
Ciach also describes $t_m$ by the translation-invariant metric
$d_m(x,y)=\rho_m(x-y)$, where
\[
 \rho_m(z):=\inf\{r>0:m(q_r(z))\le r\}
\]
for $z\in L_m(R)$; see \cite[Definition~2.5 and
Proposition~2.6]{Ciach}. Theorem~\ref{thm:topology} asserts
\emph{equality of topologies}, not equality of the numerical metrics
$\rho_m$, $\rho_\tau$, or the operator norm.
\end{rem}

\section{Completeness, density and strictness}
\begin{thm}[Completeness and completion of the bounded part]
\label{thm:completion}
For every Ciach measure $m$, the algebra $(LL_m(R),t_m)$ is a
complete metrizable topological $*$-algebra. Moreover:
\begin{enumerate}\renewcommand{\labelenumi}{(\roman{enumi})}
\item If $\delta_m>0$, this is the Banach $*$-algebra $R$ with its
operator norm topology.
\item If $\delta_m=0$, this is $(S(R),t_\tau)$, and the canonical
inclusion
\[
 \iota:(R,t_{\|\cdot\|})\longrightarrow(S(R),t_\tau)
\]
is continuous, injective and has dense, proper range.
\item The relative topology $t_\tau|_R$ is strictly weaker than
$t_{\|\cdot\|}$. The space $(R,t_\tau|_R)$ is not complete,
and its Hausdorff completion as a topological vector space is
$(S(R),t_\tau)$.
\end{enumerate}
\end{thm}
\begin{proof}
Completeness, metrizability and continuity of the algebra operations
follow from Ciach's Corollaries~2.5--2.7 applied to the algebra
$L_m(R)$, together with Theorem~\ref{thm:topology}; alternatively,
for $m=\tau$, this is the standard theorem of Nelson
\cite{Nelson}. Part (i) now follows from completeness of $R$ in
operator norm.

For part (ii), note first that if $x\in R$ and
$\|x\|\le\varepsilon$, then $q_\varepsilon(x)=0$, so
$x\in V_\tau(\varepsilon,\eta)$ for every $\eta>0$.
Thus $\iota$ is continuous. To prove density, let $x\in S(R)$
and define its bounded spectral truncations
\begin{equation}\label{eq:truncations}
 x_N:=xE_{[0,N]}(|x|)\in R,\qquad N\ge1.
\end{equation}
The spectral theorem and normality of $\tau$ give
\[
 \tau\bigl(E_{(N,\infty)}(|x|)\bigr)\downarrow0.
\]
For every fixed $\varepsilon>0$ and $N>\varepsilon$,
$q_\varepsilon(x-x_N)=E_{(N,\infty)}(|x|)$.
Indeed, in the polar decomposition $x=v|x|$, the operator
$x-x_N=v|x|E_{(N,\infty)}(|x|)$, so its absolute value is
$|x|E_{(N,\infty)}(|x|)$.
By Lemma~\ref{lem:spectral} with $m=\tau$,
$x_N\to x$ in $t_\tau$. Thus $R$ is dense in $S(R)$.

For strictness, choose mutually orthogonal projections
$(e_n)_{n\ge1}\subseteq R$ with
$\tau(e_n)=2^{-n}$ and $\displaystyle\bigvee_{n\ge1}e_n=\e$.
Such a family exists by diffuseness and normality of the trace.
Then
\begin{equation}\label{eq:unbounded-h}
 h:=\sum_{n\ge1}n e_n\in S(R)\setminus R.
\end{equation}
Thus the inclusion is proper.

For part (iii), choose projections $p_n\neq0$ with
$\tau(p_n)=2^{-n}$. For every $0<\varepsilon<1$,
$q_\varepsilon(p_n)=p_n$, and therefore
\[
 p_n\xrightarrow[t_\tau]{}0,
 \qquad\|p_n\|=1\quad(n\ge1).
\]
This proves that $t_\tau|_R$ is strictly weaker than the norm
topology. The truncated sequence
$h_N=\sum_{n=1}^N n e_n\in R$ converges in $t_\tau$ to the
unbounded operator $h$ from \eqref{eq:unbounded-h}, hence is
$t_\tau|_R$-Cauchy. If it converged in $R$, the Hausdorff property
of the $t_\tau$-topology on $S(R)$ would force its limit to equal
$h\notin R$. Hence $(R,t_\tau|_R)$ is incomplete.
Finally, a dense subspace of a complete Hausdorff topological
vector space has that space as its Hausdorff completion, proving
(iii).
\end{proof}

\begin{cor}[Projection convergence]
\label{cor:projections}
Let $(p_n)$ be a sequence of projections in $R$.
Then
\[
 p_n\xrightarrow[t_m]{}0
 \quad\Longleftrightarrow\quad m(p_n)\to0.
\]
If $\delta_m=0$, this is equivalent to $\tau(p_n)\to0$.
If $\delta_m>0$, it is equivalent to $p_n=0$ for all
sufficiently large $n$.
\end{cor}
\begin{proof}
For $0<\varepsilon<1$, $q_\varepsilon(p_n)=p_n$, and
for $\varepsilon\ge1$, $q_\varepsilon(p_n)=0$.
The spectral criterion of Lemma~\ref{lem:spectral} yields the
first assertion. The remaining assertions follow from
Theorem~\ref{thm:topology} or directly from the inequalities
$m(p_n)\ge\delta_m$ for $p_n\neq0$ and the relation
$m(p_n)=f(\tau(p_n))$ when $\delta_m=0$.
\end{proof}

\subsection*{Topological rigidity of norm-valued linear maps.}

The following obstruction holds for \emph{all} continuous linear maps,
not merely for derivations. It is useful for distinguishing the
domain and range topologies correctly.

\begin{thm}[No nonzero continuous norm-valued linear maps]
\label{thm:zero-dual}
Let $Y$ be any normed complex vector space, not necessarily complete.
Let $X$ be either $S(R)$ with its $\tau$-measure topology or $R$
with the relative topology $t_\tau|_R$.
If $T:X\to Y$ is complex-linear and continuous, then
\begin{equation}\label{eq:Tzero}
 T=0.
\end{equation}
In particular, the continuous dual of $(S(R),t_\tau)$ and
that of $(R,t_\tau|_R)$ are both $\{0\}$.
\end{thm}
\begin{proof}
Continuity at zero implies that there exist $\varepsilon,\eta>0$
such that
\begin{equation}\label{eq:bounded-T}
 \|T(y)\|_Y<1\qquad
 \bigl(y\in V_\tau(\varepsilon,\eta)\cap X\bigr).
\end{equation}
Choose an integer $N$ such that $1/N<\eta$. Since $R$ is
a type $\mathrm{II}_1$ factor, there exist pairwise orthogonal
projections $p_1,\ldots,p_N\in R$ satisfying
\[
 \sum_{j=1}^Np_j=\e,\qquad \tau(p_j)=1/N<\eta.
\]
Fix any $x\in X$ and $j\in\{1,\ldots,N\}$. The operator
$x p_j$ belongs to $X$ (either $X=S(R)$ or $X=R$).
For every $\lambda\in\mathbb{C}$, putting $y=\lambda x p_j$ in the
algebra $S(R)$, we have $y=yp_j$. The right support of $y$
(the support of $|y|$) is therefore dominated by $p_j$:
this follows from $y^*y=p_jy^*yp_j$, interpreted as an equality
of affiliated operators. In particular, this remains valid if $x$
is unbounded. Therefore, for every $\varepsilon>0$,
\[
 q_\varepsilon(\lambda x p_j)\le p_j,
 \qquad
 \tau\bigl(q_\varepsilon(\lambda x p_j)\bigr)
 \le\tau(p_j)<\eta.
\]
By Lemma~\ref{lem:spectral},
$\lambda x p_j\in V_\tau(\varepsilon,\eta)\cap X$
for \emph{every} $\lambda\in\mathbb{C}$. Thus \eqref{eq:bounded-T}
gives
\[
 |\lambda|\,\|T(x p_j)\|_Y<1
 \qquad(\lambda\in C).
\]
Letting $|\lambda|\to\infty$ shows $T(x p_j)=0$.
Since $\sum_jp_j=1$, multiplication in the affiliated operator
algebra yields $x=\sum_jx p_j$, and therefore
\[
 T(x)=\sum_{j=1}^NT(x p_j)=0.
\]
As $x\in X$ was arbitrary, $T=0$. Taking $Y=C$ gives the
continuous-dual assertion.
\end{proof}

\begin{cor}[Non-local-convexity]
\label{cor:nonloc}
Neither $(S(R),t_\tau)$ nor $(R,t_\tau|_R)$ is a locally convex
topological vector space. In particular, the two regimes in
Theorem~\ref{thm:topology} cannot be identified by a topological
vector-space isomorphism.
\end{cor}
\begin{proof}
Both spaces are nonzero and Hausdorff. By
Theorem~\ref{thm:zero-dual}, they have trivial continuous dual.
A nonzero Hausdorff locally convex complex vector space admits
nonzero continuous linear functionals separating points (by the
Hahn--Banach separation theorem). Hence neither space can be
locally convex. By contrast, $(R,t_{\|\cdot\|})$ is a Banach space
and therefore locally convex. Local convexity is preserved by
topological vector-space isomorphisms.
\end{proof}

\subsection*{Automatic continuity of derivations.}

By a derivation we always mean a complex-linear map satisfying
the Leibniz identity. No boundedness or continuity is assumed
a priori.

\begin{thm}[Automatic $m$-measure continuity]
\label{thm:automatic}
For every subadditive measure $m$ on a type
$\mathrm{II}_1$ factor $R$, every derivation
\[
 D:LL_m(R)\longrightarrow LL_m(R)
\]
is continuous with respect to $t_m$. More precisely:

\begin{enumerate}
\item If $\delta_m>0$, there exists $a\in R$ with
$D(x)=[a,x]$, and $D$ is operator-norm bounded, with
\[
 \|D(x)\|\le2\|a\|\,\|x\|\qquad(x\in R).
\]
\item If $\delta_m=0$, there exists $a\in S(R)$ with
$D(x)=[a,x]$ for $x\in S(R)$, and $D$ is continuous in
the measure topology $t_\tau$.
\end{enumerate}
Moreover, every complex-linear derivation
$D:R\to S(R)$ is continuous as a map
$(R,t_\tau|_R)\to(S(R),t_\tau)$ and, consequently, as a map
$(R,t_{\|\cdot\|})\to(S(R),t_\tau)$.
\end{thm}
\begin{proof}
If $\delta_m>0$, Theorem~\ref{thm:topology} identifies the
domain and codomain with $R$ in operator norm. Sakai's theorem
\cite{Sakai} provides $a\in R$ with $D=ad_a$, and
\[
 \|[a,x]\|\le\|ax\|+\|xa\|
 \le2\|a\|\,\|x\|.
\]
This proves (i).

If $\delta_m=0$, Theorem~\ref{thm:topology} identifies the
algebra with $(S(R),t_\tau)$. The theorem of
Ber--Kudaybergenov--Sukochev \cite{BKS} gives
$D=ad_a$ for some $a\in S(R)$. Ciach's multiplication
continuity theorem \cite[Corollary~2.6]{Ciach} guarantees
that the maps $x\mapsto ax$ and $x\mapsto xa$ are continuous
in $t_\tau$, hence so is $D=ad_a$. This proves (ii).
Finally, by Huang--Kudaybergenov--Sukochev
\cite[Theorem~1.1]{HKS}, every derivation $D:R\to S(R)$ has
the form $ad_a|_R$ for some $a\in S(R)$. Since
$ad_a:S(R)\to S(R)$ is $t_\tau$-continuous, its restriction
to $R$ is continuous in the relative topology. The norm-to-
measure continuity follows from the continuity of the inclusion
established in Theorem~\ref{thm:completion}.
\end{proof}

\begin{cor}[A sharp domain/range topology distinction]
\label{cor:range-norm}
Let $D:R\to R$ be a complex-linear derivation. Then
\[
 D:(R,t_\tau|_R)\longrightarrow(R,t_{\|\cdot\|})
 \text{ is continuous}
 \quad\Longleftrightarrow\quad D=0.
\]
Nevertheless every such $D$ is continuous both as
$(R,t_{\|\cdot\|})\to(R,t_{\|\cdot\|})$ and as
$(R,t_\tau|_R)\to(R,t_\tau|_R)$.
\end{cor}
\begin{proof}
The first assertion follows immediately from
Theorem~\ref{thm:zero-dual}, taking $Y=(R,\|\cdot\|)$.
The converse is trivial. The norm continuity follows from
Sakai's theorem and the estimate in
Theorem~\ref{thm:automatic}(i); the measure continuity follows
from multiplication continuity, since a norm-implementer
$a\in R$ also belongs to $S(R)$.
\end{proof}

\begin{prop}[An explicit discontinuity witness]
\label{prop:explicit}
Let $a\in R\setminus \mathbb{C}\e$ and $D=ad_a|_R$. Then there
exists a sequence of nonzero projections $e_n\in P(R)$ such
that
\[
 \tau(e_n)\longrightarrow0,
 \qquad
 \inf_{n\ge1}\|D(e_n)\|>0.
\]
Thus $e_n\to0$ in $t_\tau$, while $D(e_n)$ does not converge
to zero in operator norm.
\end{prop}
\begin{proof}
Because $a$ is not central, it fails to commute with some
projection $p\in P(R)$. Replacing $p$ by $\e-p$, if needed,
we may arrange that $b:=pa(\e-p)\neq0$.
Choose $c>0$ such that the spectral projection
\[
 e:=E_{[c,\infty)}(|b|)\neq0.
\]
Since $b=b(\e-p)$, one has $e\le\e-p$. By diffuseness of the
corner $eRe$, choose nonzero projections $e_n\le e$ with
$\tau(e_n)\to0$ (for instance, $\tau(e_n)=2^{-n}\tau(e)$).
For every $n$, the bound $|b|e\ge c e$ implies
\begin{equation}\label{eq:be-lower}
 \|b e_n\|\ge c.
\end{equation}
Indeed, for each unit vector $\xi\in e_nH$,
$\|b\xi\|=\||b|\xi\|\ge c$.
Since $e_n\le\e-p$, we have $pe_n=0$ and
\[
 p[a,e_n]e_n=pae_n=b e_n.
\]
Therefore, by \eqref{eq:be-lower},
\[
 \|[a,e_n]\|\ge\|p[a,e_n]e_n\|
 =\|b e_n\|\ge c>0.
\]
On the other hand, for $0<\varepsilon<1$,
$\tau(q_\varepsilon(e_n))=\tau(e_n)\to0$.
Hence $e_n\to0$ in $t_\tau$ by Lemma~\ref{lem:spectral}.
\end{proof}

%\section{Examples, comparison and scope}

\begin{ex}[Four explicit subadditive measures]
For $0<\alpha\le1$, consider
\[
 \begin{array}{ll}
 m_\tau(p)=\tau(p),&
 m_\alpha(p)=\tau(p)^\alpha,\\[2pt]
 m_0(p)=\e_{\{p\neq0\}},&
 m_{0,\alpha}(p)=\e_{\{p\neq0\}}+\tau(p)^\alpha.
 \end{array}
\]
These measures satisfy Ciach's axioms. Their small-projection
gaps are, respectively, $0,0,1,1$. Consequently,
\[
 \begin{aligned}
 (LL_{m_\tau}(R),t_{m_\tau})
 &=(LL_{m_\alpha}(R),t_{m_\alpha})=(S(R),t_\tau),\\
 (LL_{m_0}(R),t_{m_0})
 &=(LL_{m_{0,\alpha}}(R),t_{m_{0,\alpha}})
 =(R,t_{\|\cdot\|}).
 \end{aligned}
\]
The equality of topologies does not assert equality of their
corresponding numerical measure metrics.
\end{ex}

\begin{cor}[Complete topology/derivation dichotomy]
Let $m$ be an arbitrary Ciach subadditive measure on a type
$\mathrm{II}_1$ factor. Then
\[
 \boxed{
 \begin{array}{c|c|c|c}
 &LL_m(R)&t_m& Der(LL_m(R))\\ \hline
 \delta_m>0&R&\text{operator norm}&
 \{ad_a:a\in R\}\\[2pt]
 \delta_m=0&S(R)&\tau\text{-measure}&
 \{ad_a:a\in S(R)\}
 \end{array}}
\]
and in both rows all derivations are automatically $t_m$-continuous.
In the second row the continuous dual of the underlying topological
vector space is zero; in the first row the algebra is a Banach space.
\end{cor}
\begin{proof}
Combine Theorems~\ref{thm:topology}, \ref{thm:zero-dual}
and \ref{thm:automatic}, together with Sakai's and the
Ber--Kudaybergenov--Sukochev innerness theorems.
\end{proof}

\begin{rem}[Established inputs and scope]
The $m$-measure topology, its completeness, its metrizability
and its topological $*$-algebra structure are classical results
of Ciach \cite[Section~2, Proposition~2.6 and Corollaries~2.5--2.7]{Ciach}.
The spectral-neighborhood equivalence proved here for a type
$\mathrm{II}_1$ factor is also consistent with Ciach's general
formula for generalized singular numbers
\cite[Proposition~2.3]{Ciach}. The usual measure topology was studied by
Nelson \cite{Nelson}; the innerness assertions used here are due
to Sakai \cite{Sakai}, Ber--Kudaybergenov--Sukochev
\cite{BKS}, and Huang--Kudaybergenov--Sukochev \cite{HKS}.
The two-regime conclusion combines those known results with the
Stage II algebraic dichotomy. The vanishing continuous-dual
argument and the explicit discontinuity witness have been included
with complete proofs, but their inclusion is not a claim of
historical novelty. For an original result beyond these factor
reductions, the nonfactorial type II setting, where central
localization may produce intermediate algebras, is a more
promising setting.
\end{rem}

\section{Central localization and derivations for sub\-additive-measure operator algebras.}

Let $M\subseteq B(H)$ be a complex von Neumann algebra. Define
\[
 \begin{aligned}
E(M)=\bigl\{x\in LS(M):\;&\exists\text{ pairwise orthogonal }
 (e_j)_{j\in J}\subseteq P(Z(M)),\\
 &\bigvee_{j\in J}e_j=\e,\quad e_jx\in M\ (j\in J)\bigr\}.
 \end{aligned}
\]
The index set $J$ need not be countable. For clarity, a derivation
is a complex-linear map satisfying the Leibniz identity, unless an
\emph{additive ring derivation} is explicitly mentioned.
%\section{Universal central-localization statements}

\begin{lem}[Countable central localization is transitive]
\label{lem:trans}
Suppose $x$ is affiliated with $M$ and $z_n\uparrow\e$ are central
projections with $z_nx\in LS(M)$ for all $n$. Then $x\in LS(M)$.
\end{lem}
\begin{proof}
For every $n$ choose central projections $w_{n,k}\uparrow\e$ such
that $w_{n,k}z_nx\in S(M)$. Set $z_0=0$ and $e_n=z_n-z_{n-1}$.
For $N\in N$, put
\[
 u_N=\sum_{n=1}^{N}e_nw_{n,N}.
\]
The summands are orthogonal central projections and
$u_N\le u_{N+1}$. For fixed $n$,
$\bigvee_{N\ge n}e_nw_{n,N}=e_n$, hence $u_N\uparrow\e$.
Moreover,
\[
 u_Nx=\sum_{n=1}^{N}e_nw_{n,N}z_nx\in S(M),
\]
using stability of $S(M)$ under central cuts and finite sums.
Thus the sequence $(u_N)$ witnesses $x\in LS(M)$.
\end{proof}

\begin{lem}[The $m$-measurable algebra is locally measurable]
\label{lem:LintoLS}
For every von Neumann algebra $M$ and every Ciach measure $m$,
\[
 L_m(M)\subseteq LS(M).
\]
\end{lem}
\begin{proof}
We recall the standard spectral obstruction to local
measurability, in the form of Mori's criterion
\cite[Corollary~2.2]{Mori}. If a positive affiliated operator
$y$ does not belong to $LS(M)$, there exist numbers
$0<c_1<c_2<\cdots$ tending to infinity and nonzero,
mutually Murray--von Neumann equivalent projections
$f_n\le E_{[c_n,c_{n+1})}(y)$. If $x\in L_m(M)$ but
$y=|x|\notin LS(M)$, then for every $n$,
\[
 0<m(f_1)=m(f_n)\le m\bigl(E_{(c_n/2,\infty)}(y)\bigr),
\]
contradicting the spectral-tail criterion for $x\in L_m(M)$. Thus $|x|\in LS(M)$, and the
polar decomposition $x=v|x|$ with $v\in M$ implies
$x\in LS(M)$. This argument also covers $m(f_1)=\infty$.
\end{proof}

\begin{thm}[Universal sandwich and the indicator measure]
\label{thm:sandwich}
For every von Neumann algebra $M$ and every Ciach measure $m$,
\begin{equation}\label{eq:sandwich}
E(M)\subseteq LL_m(M)\subseteq LS(M).
\end{equation}
Let $m_0(p)=0$ for $p=0$ and $m_0(p)=1$ otherwise. Then $m_0$
is a Ciach measure and
\begin{equation}\label{eq:indicator}
L_{m_0}(M)=M,\qquad LL_{m_0}(M)=E(M).
\end{equation}
No $\sigma$-finiteness assumption is required.
\end{thm}
\begin{proof}
For $x\in LL_m(M)$, take $z_n\uparrow1$ with
$z_nx\in L_m(M)$. Lemmas~\ref{lem:LintoLS} and~\ref{lem:trans}
imply $x\in LS(M)$.

Now take $x\in E(M)$, witnessed by a possibly uncountable central
partition $(e_j)_{j\in J}$. Set
$\alpha_j=\|e_jx\|<\infty$ and
\[
 z_n=\bigvee_{\{j:\alpha_j\le n\}}e_j.
\]
Then $z_n\uparrow\e$. As $|z_nx|\le nz_n$, the operator $z_nx$ is
bounded and affiliated, hence belongs to $M\subseteq L_m(M)$.
This proves the left inclusion in \eqref{eq:sandwich}.

For $m_0$, faithfulness, monotonicity, equivalence invariance,
subadditivity, and sequential normality are immediate. By the spectral-tail criterion for $m$-measurability, we obtain $x\in L_{m_0}(M)$ if and only if its
spectral tail vanishes for some finite threshold, which is equivalent
to $x$ being bounded. Thus $L_{m_0}(M)=M$.

If $x\in LL_{m_0}(M)$, write $z_n\uparrow\e$ with $z_nx\in M$.
The disjoint central projections $e_1=z_1$ and
$e_n=z_n-z_{n-1}$ for $n\ge2$ form a partition of unity with
$e_nx\in M$. Hence $x\in E(M)$. The reverse inclusion is the
first part of the theorem, giving \eqref{eq:indicator}.
\end{proof}

\begin{prop}[Type II forces a strict upper inclusion for $m_0$]
\label{prop:strictE}
If $M$ has a nonzero type II direct summand, then
\[
E(M)=LL_{m_0}(M)\subsetneq LS(M).
\]
Moreover one can choose the witness in $S(M)\setminus E(M)$.
\end{prop}
\begin{proof}
Choose a nonzero finite projection $p$ inside the type II
summand. By repeated halving in the continuous finite corner
$pMp$, obtain mutually orthogonal nonzero projections $e_n\le p$
with common, nonzero central support $c_M(e_n)=c_M(p)$.
Explicitly, put $r_0=p$ and recursively split
$r_{n-1}=e_n+r_n$ into equivalent orthogonal projections, so that
$c_M(e_n)=c_M(r_n)=c_M(p)$.
Set
\[
 h=\sum_{n\ge1}n e_n,
\]
with value zero on the orthogonal complement. It is a densely
defined positive self-adjoint affiliated operator. Since
$E_{(N,\infty)}(h)=\bigvee_{n>N}e_n\le p$, its high spectral
tails are finite; hence $h\in S(M)$.

If $h\in E(M)$, some nonzero central cut $zh$ with
$zc_M(p)\ne0$ would have to be bounded. But
$ze_n\ne0$ for every $n$ by equality of central supports, and
$zh$ has the eigenvalue $n$ on $ze_nH$. Thus
$\|zh\|\ge n$ for every $n$, a contradiction. Consequently
$h\notin E(M)$.
\end{proof}
%\section{The two-regime theorem for both kinds of type II factors}
\begin{lem}[Trace and Segal measurability in a type II factor]
\label{lem:factortrace}
Let $F$ be a type $\mathrm{II}_1$ or $\mathrm{II}_{\infty}$
factor with a faithful normal semifinite trace $\tau$ (normalized
when $F$ is finite). For a densely defined affiliated operator $x$,
\[
 x\in S(F)
 \quad\Longleftrightarrow\quad
 \exists t_0>0:\ \tau\bigl(E_{(t_0,\infty)}(|x|)\bigr)<\infty.
\]
If $x\in S(F)$, then
$\tau\bigl(E_{(t,\infty)}(|x|)\bigr)\longrightarrow0$ \ ($t\to\infty$).
\end{lem}
\begin{proof}
In a semifinite factor, a projection is Murray--von Neumann
finite if and only if its trace is finite. If
$q_{t_0}:=E_{(t_0,\infty)}(|x|)$ has finite trace, then
$E_{[0,t]}(|x|)\uparrow\e$, its ranges lie in $\mathcal D(x)$,
and its complements are finite for every $t\ge t_0$.
Thus $x\in S(F)$.

Conversely, suppose $x\in S(F)$. By Segal measurability,
there are projections $p_n\uparrow\e$ with
$p_nH\subseteq\mathcal D(x)$ and $\e-p_n$ finite. Since $x$
is closed and $xp_n$ is defined on all of $H$, the closed graph
theorem gives $\|xp_n\|<\infty$. Fix $n$ and choose
$t>\|xp_n\|$. Then
\[
 p_n\wedge E_{(t,\infty)}(|x|)=0.
\]
The projection-lattice equivalence
$r-(r\wedge s)\sim_F (r\vee s)-s$ implies that
$E_{(t,\infty)}(|x|)\precsim_F \e-p_n$.
Hence $q_t$ is finite and has finite trace.
Finally $q_s\downarrow0$ strongly as $s\to\infty$, and
$q_s\le q_t$ for $s\ge t$. Normality of the trace, applied
inside the finite-trace corner $q_tFq_t$, yields
$\tau(q_s)\downarrow0$.
\end{proof}

\begin{thm}[Type $\mathrm{II}_{1}$ and type $\mathrm{II}_{\infty}$ factor dichotomy]
\label{thm:factor}
Let $F$ be a type $\mathrm{II}_1$ or type $\mathrm{II}_{\infty}$
factor. For any Ciach measure $m$ on $P(F)$, define
\[
 \delta_m=\inf\{m(p):0\ne p\in P(F)\}\in[0,\infty].
\]
Then
\begin{equation}\label{eq:factordichotomy}
 LL_m(F)=L_m(F)=
 \begin{cases}
 F,&\delta_m>0,\\[2pt]
 S(F)=LS(F),&\delta_m=0.
 \end{cases}
\end{equation}
The result holds without a separability assumption, and
$\delta_m=\infty$ is included in the first case.
\end{thm}
\begin{proof}
Since $F$ is a factor, the only central projections are $0$ and
$1$; hence central localization does not change the algebra:
$LL_m(F)=L_m(F)$. Also $LS(F)=S(F)$.

If $\delta_m>0$ and $x\in L_m(F)$, choose a finite real
number $\varepsilon$ with $0<\varepsilon<\delta_m$ (also
possible for $\delta_m=\infty$). The spectral-tail criterion for $x\in L_m(F)$ forces
$m(E_{(t,\infty)}(|x|))<\varepsilon<\delta_m$
for sufficiently large $t$.
Therefore that spectral projection is zero, so $x$ is bounded.
Conversely bounded affiliated operators belong to $L_m(F)$.
Thus $L_m(F)=F$.

Suppose $\delta_m=0$. Let $x\in S(F)$ and $\varepsilon>0$.
There is a nonzero projection $r\in P(F)$ such that
$0<m(r)<\varepsilon$. By semifiniteness choose a nonzero
\emph{finite} subprojection $e\le r$; monotonicity gives
$0<m(e)<\varepsilon$, and $0<\tau(e)<\infty$.
By Lemma~\ref{lem:factortrace}, for large $t$ we have
\[
 \tau\bigl(q_t\bigr)\le\tau(e),\qquad
 q_t=E_{(t,\infty)}(|x|).
\]
Projection comparison by trace in a semifinite factor gives
$q_t\precsim_F e$. Hence $q_t\sim_F e'_t$ for some
$e'_t\le e$. Invariance and monotonicity of $m$ yield
\[
 m(q_t)=m(e'_t)\le m(e)<\varepsilon.
\]
As $\varepsilon$ is arbitrary, the spectral-tail criterion implies
$x\in L_m(F)$. Thus $S(F)\subseteq L_m(F)$.
The reverse inclusion follows from Lemma~\ref{lem:LintoLS}
and $LS(F)=S(F)$. This proves the theorem.
\end{proof}

\begin{rem} This Theorem is proved in Theorem \ref{thm:dichotomy} by another method for type II$_1$.
Unlike the type $\mathrm{II}_1$ case, $S(F)$ for a type
$\mathrm{II}_{\infty}$ factor is usually a proper subset of all
affiliated operators. Formula \eqref{eq:factordichotomy} does
\emph{not} assert $S(F)=Aff(F)$ in the infinite case.
\end{rem}

\subsection*{Classification for a countable atomic center.}

Let $I$ be finite or countably infinite, let $F_i$ be type II
factors (each of type $\mathrm{II}_1$ or $\mathrm{II}_{\infty}$), and
let $\displaystyle M=\prod_{i\in I}^{\ell^\infty} F_i$ be their von Neumann
algebraic direct product. The central atoms $z_i$ are the coordinate
units. We use $\prod$ without the $\ell^\infty$ superscript for
\emph{unrestricted} products of affiliated-operator algebras.

\begin{lem}[Coordinatewise locally measurable operators]
\label{lem:LSproduct}
For a finite or countable product of type II factors,
\[
 E(M)=\prod_{i\in I}F_i,
 \qquad
 LS(M)=\prod_{i\in I}S(F_i).
\]
These are unrestricted algebraic products on the right.
\end{lem}
\begin{proof}
For $x\in LS(M)$, taking the central cut $z_i$ shows that its
$i$th coordinate is in $S(F_i)$; equivalently, a suitable central
localization sequence must eventually contain the central atom
$z_i$. Thus $\displaystyle LS(M)\subseteq\prod_iS(F_i)$.
Conversely, for $x=(x_i)$ with $x_i\in S(F_i)$, enumerate
$I$ and let $\displaystyle w_N=\sum_{i\le N}z_i$. Then $w_N\uparrow\e$, and
$w_Nx$ is a finite direct sum of measurable operators, hence
belongs to $S(M)$. Thus $x\in LS(M)$.

Similarly, $x\in E(M)$ implies every $z_ix$ is bounded; each
coordinate belongs to $F_i$. Conversely, for $(a_i)$ with
$a_i\in F_i$, group the indices according to $\|a_i\|\le N$.
Each corresponding central cut is bounded, and these increasing
cuts tend to $1$. By the definition of $E(M)$ (or the bounded-cut
criterion in Theorem~\ref{thm:sandwich}), $x\in E(M)$.
\end{proof}

\begin{lem}[Restriction of an arbitrary Ciach measure]
\label{lem:restr}
Let $m$ be a Ciach measure on $P(M)$, without any additive
or product-decomposition hypothesis. Define
\[
 m_i(p)=m(0,\ldots,0,p,0,\ldots),\qquad p\in P(F_i).
\]
Then every $m_i$ is a Ciach measure on $F_i$.
\end{lem}
\begin{proof}
Each of Ciach's five axioms passes to the corner $z_iM$ under
its projection-lattice identification with $F_i$. In particular,
Murray--von Neumann equivalence and increasing suprema of
projections in the corner agree with those in $M$.
\end{proof}

\begin{thm}[Exact coordinatewise classification]
\label{thm:countable}
With the assumptions above, put
\[
 \delta_i=\inf\{m_i(p):0\ne p\in P(F_i)\},\qquad
 A_i=\begin{cases}
 F_i,&\delta_i>0,\\
 S(F_i),&\delta_i=0.
 \end{cases}
\]
Then for \emph{every} Ciach measure $m$ on $M$,
\begin{equation}\label{eq:product}
LL_m(M)=\prod_{i\in I}A_i.
\end{equation}
Neither finiteness of $m(\e)$ nor separability of the $F_i$
is required. In particular, $LL_m(M)$ depends only on which
restricted thresholds $\delta_i$ vanish.
\end{thm}
\begin{proof}
Let $x\in LL_m(M)$ and choose central projections
$w_N\uparrow\e$ with $w_Nx\in L_m(M)$. For each fixed $i$,
$z_i\le w_N$ for all sufficiently large $N$, since the $i$th
coordinate of an increasing sequence of central projections
can take only the values $0$ and $1$. Consequently
$z_ix\in L_m(M)$ by stability under central cuts. Its spectral
tails lie in the $i$th corner. Applying the spectral-tail criterion to $z_ix$ gives
\[
 \lim_{t\to\infty}m_i(E_{(t,\infty)}(|x_i|))=0,
\]
and hence $x_i\in L_{m_i}(F_i)=A_i$ by
Theorem~\ref{thm:factor}. This proves one inclusion.

Conversely, suppose $x=(x_i)$ with every $x_i\in A_i$.
Enumerate $I$ and put $w_N=\sum_{i\le N}z_i$. For each $N$,
the spectral tail of $w_Nx$ is the join of the finitely many
coordinate spectral tails. By subadditivity,
\[
 \begin{aligned}
 m\bigl(E_{(t,\infty)}(|w_Nx|)\bigr)\le\sum_{i\le N}m_i\bigl(E_{(t,\infty)}(|x_i|)\bigr)
 \longrightarrow0\quad(t\to\infty).
 \end{aligned}
\]
Thus $w_Nx\in L_m(M)$ by the spectral-tail criterion~\cite[Proposition~2.1]{Ciach}; since
$w_N\uparrow\e$, we conclude $x\in LL_m(M)$. The argument
uses only finite subadditivity, not any decomposition formula for $m$.
\end{proof}

\begin{cor}[Center and derivations of the atomic product]
\label{cor:atomicder}
In the setting of Theorem~\ref{thm:countable},
\begin{equation}\label{eq:centerprod}
 Z(LL_m(M))=\prod_{i\in I} \mathbb{C}\e_i
 \cong C^{I}
\end{equation}
and
\begin{equation}\label{eq:derprod}
 Der(LL_m(M))\cong_{\mathrm{Lie}}
 \prod_{i\in I}(A_i/ \mathbb{C}\e_i).
\end{equation}
Every derivation is inner, implemented by an element of
$LL_m(M)$; no continuity hypothesis is imposed.
\end{cor}
\begin{proof}
By \eqref{eq:product}, the center is the unrestricted product of
the centers of $A_i$. Since $A_i$ is either $F_i$ or $S(F_i)$
and $F_i$ is a factor, $Z(A_i)=\mathbb{C}\e_i$, giving
\eqref{eq:centerprod}.

Let $D$ be a derivation of $\displaystyle A=\prod_iA_i$. Each central
idempotent $z_i$ satisfies $D(z_i)=0$: in a unital algebra a
derivation sends central elements into the center, and
$D(z_i)=D(z_i^2)=2z_iD(z_i)$ implies $D(z_i)=0$.
Therefore $D(z_ix)=z_iD(x)$, and $D$ is coordinatewise:
$D(x)_i=D_i(x_i)$ for a derivation $D_i:A_i\to A_i$.

If $A_i=F_i$, Sakai's theorem makes $D_i$ inner. If
$A_i=S(F_i)$, the known innerness theorem for type
$\mathrm{II}_1$ Murray--von Neumann algebras, or the theorem
for properly infinite algebras in the type $\mathrm{II}_{\infty}$
case, makes $D_i$ inner \cite{Sakai,BKS2,BChS}.
Choose $a_i\in A_i$ with $D_i=ad_{a_i}$.
The unrestricted family $a=(a_i)$ belongs to $A$, and
$D=ad_a$. Its kernel as a function of $a$ is precisely
$Z(A)$. The first isomorphism theorem for complex Lie
algebras yields \eqref{eq:derprod}.
\end{proof}
%\section{A strict intermediate algebra already on $R\oplus R$}

\begin{thm}[The minimal nonfactor example]
\label{thm:two}
Let $R$ be a type $\mathrm{II}_1$ factor with normalized trace
$\tau$, and put $M=R\oplus R$. Define
\[
 m(p\oplus q)=\tau(p)+\e_{\{q\ne0\}}.
\]
Then $m$ is a Ciach measure and
\begin{equation}\label{eq:intermediate}
 E(M)=R\oplus R
 \ \subsetneq\
 LL_m(M)=S(R)\oplus R
 \ \subsetneq\
 LS(M)=S(R)\oplus S(R).
\end{equation}
Moreover
\begin{equation}\label{eq:twoder}
 Der(LL_m(M))\cong_{\mathrm{Lie}}
 (S(R)/ \mathbb{C}\e)\oplus(R/ \mathbb{C}\e).
\end{equation}
\end{thm}
\begin{proof}
The trace term and indicator term separately satisfy
faithfulness on their nonzero corners, monotonicity,
invariance under equivalence, subadditivity and sequential
normality. Their sum is faithful on $M$ and satisfies all
five axioms. The coordinate restrictions have thresholds
$\delta_1=0$ and $\delta_2=1$. Thus
Theorem~\ref{thm:countable} gives
$LL_m(M)=S(R)\oplus R$, while Lemma~\ref{lem:LSproduct}
gives the two outer algebras in \eqref{eq:intermediate}.

Choose an unbounded positive operator $h\in S(R)\setminus R$,
for example $h=\sum_{n\ge1}ne_n$ with orthogonal projections
$\tau(e_n)=2^{-n}$. Then $h\oplus0$ lies in
$LL_m(M)\setminus E(M)$, whereas $0\oplus h$ lies in
$LS(M)\setminus LL_m(M)$. Both inclusions are strict.
Formula \eqref{eq:twoder} follows from
Corollary~\ref{cor:atomicder}.
\end{proof}

\begin{rem}
The example proves that a vanishing global number
$\inf_{p\ne0}m(p)$ does \emph{not} determine $LL_m(M)$ when
the center is nontrivial: here the global infimum is zero,
yet $LL_m(M)\ne LS(M)$. The separate central-corner
thresholds, rather than a single global threshold, are needed.
\end{rem}

%\section{Why the countability hypothesis cannot be suppressed}

\begin{thm}[Uncountable atomic-center obstruction]
\label{thm:uncountable}
Let $R$ be a type $\mathrm{II}_1$ factor, set
$I=N^{N}$, and let
\[
 M=\prod_{f\in I}^{\ell^\infty}R,
 \qquad
 m(p)=\sup_{f\in I}\tau(p_f),
 \quad p=(p_f)\in P(M).
\]
Then $m$ is a finite, faithful Ciach measure, and its
restriction to every atomic factor is the normalized trace
$\tau$. Nevertheless,
\[
LL_m(M)\subsetneq LS(M)
 =\prod_{f\in I}S(R).
\]
In particular, the coordinatewise formula
$LL_m(M)=\prod_{f\in I}L_{m_f}(R)$ is false for an
uncountable atomic center.
\end{thm}
\begin{proof}
All Ciach axioms are immediate except possibly sequential
normality. If $p_n\uparrow p$ in $M$, then
$\tau((p_n)_f)\uparrow\tau(p_f)$ for every $f$, and
\[
 \sup_n\sup_{f\in I}\tau((p_n)_f)
 =\sup_{f\in I}\sup_n\tau((p_n)_f)
 =\sup_{f\in I}\tau(p_f),
\]
which is exactly $m(p_n)\uparrow m(p)$. Also
$m(p)\le1$ for every $p$, and the restriction to an atom is
$\tau$.

Choose mutually orthogonal projections $(e_k)_{k\ge1}$ in $R$
with $\tau(e_k)=2^{-k}$. For each $f\in I=N^{N}$ define
\[
 x_f=\sum_{k\ge1} f(k)e_k\in S(R),
 \qquad x=(x_f)_{f\in I}.
\]
These spectral series define positive densely defined
self-adjoint affiliated operators. Realize the product algebra on
the Hilbert direct sum $\displaystyle H=\bigoplus_{f\in I}H_f$, and define
$\displaystyle x=\bigoplus_{f\in I}x_f$. Its domain contains the vectors with
finite coordinate support and components in $\mathcal D(x_f)$;
these form a dense subspace of $H$. Hence $x$ is a closed,
densely defined operator affiliated with $M$. The algebra $M$
is finite in the Murray--von Neumann sense: equivalence and
proper subprojection are tested coordinatewise, and every $R$
is finite. Thus every densely defined affiliated operator is
Segal-measurable (use its bounded spectral cuts, whose complements
are all finite). Consequently $x\in S(M)\subseteq LS(M)$.

Suppose $x\in LL_m(M)$. Then there are increasing central
projections $z_n\uparrow\e$ such that $z_nx\in L_m(M)$.
Since $Z(M)=\ell^\infty(I)$, we may write
$z_n=\e_{I_n}$ for sets
$I_1\subseteq I_2\subseteq\cdots$ satisfying
$\bigcup_n I_n=I$.
For every $n$ and every fixed $k$, we claim that
\begin{equation}\label{eq:boundedcoord}
 B_{n,k}:=\sup_{f\in I_n} f(k)<\infty
\end{equation}
(with supremum zero if $I_n$ is empty). Indeed, if the
supremum were infinite, then for every $t>0$ one could
choose $f\in I_n$ with $f(k)>t$. The $f$th coordinate of
$E_{(t,\infty)}(|z_nx|)$ would then dominate $e_k$, so
\[
 m\bigl(E_{(t,\infty)}(|z_nx|)\bigr)
 \ge \tau(e_k)=2^{-k}
 \qquad\text{for every }t>0.
\]
This contradicts $z_nx\in L_m(M)$ by the spectral-tail criterion.
Thus \eqref{eq:boundedcoord} holds.

Define $g\in N^{N}$ by choosing
$g(n)>B_{n,n}$ for each $n$. Since
$\bigcup_n I_n=I$, there is an $N$ with $g\in I_N$.
But then $g(N)\le B_{N,N}$ by the definition of
$B_{N,N}$, contradicting $g(N)>B_{N,N}$.
Hence $x\notin LL_m(M)$.

Every coordinate $x_f$ is in $L_{m_f}(R)=S(R)$ by the
factor dichotomy, so the asserted failure of the fiberwise
formula follows.
\end{proof}

\begin{rem}
This example distinguishes \emph{coordinatewise} small
spectral tails from the \emph{uniform} small-tail control
needed on each of countably many central localization
pieces. Even a normal, finite Ciach measure does not bridge
the gap for an uncountable atomic center. No statement about
arbitrary non-atomic direct-integral disintegrations is being
inferred from the countable product theorem.
\end{rem}

\subsection*{Derivations for arbitrary type II von Neumann algebras.}

The results in this section use two substantial \emph{external}
theorems. First, Huang--Kudaybergenov--Sukochev
\cite[Theorem~1.2]{HKS} prove that every additive ring derivation
from an $EW^*$-algebra $A$ into $LS(A_b)$ is spatially
implemented by an element of $LS(A_b)$ when the bounded part
$A_b$ has no nonzero finite type I direct summand. Second,
Ber--Sukochev's commutator estimate \cite{BS} says that for
self-adjoint $a\in LS(M)$ a self-adjoint
$c\in Z(LS(M))$ can be chosen so that for each
$\varepsilon>0$ there is a unitary $u\in M$ with
\begin{equation}\label{eq:estimate}
 |[a,u]|\ge(1-\varepsilon)|a-c|.
\end{equation}
We also use the standard characterization of $EW^*$-algebras
as solid $*$-subalgebras of $LS(M)$ containing $M$;
see \cite[Section~5]{BChS}.

\begin{lem}[Algebraic structure and the center]
\label{lem:ew}
Let $M$ be any von Neumann algebra and $m$ any Ciach measure.
Then $A=LL_m(M)$ is a unital $*$-subalgebra of $LS(M)$
containing $M$. Its bounded part is $A_b=M$; in particular,
$A$ is an $EW^*$-algebra. Moreover,
\begin{equation}\label{eq:center}
Z(A)=Z(LS(M))=S(Z(M)).
\end{equation}
\end{lem}
\begin{proof}
For central-localization sequences $(z_n)$ and $(w_n)$
witnessing $x,y\in A$, the increasing central sequence
$u_n=z_nw_n\uparrow1$ simultaneously localizes both.
Because $L_m(M)$ is a $*$-algebra stable under central cuts,
$u_n(x+y)$, $u_nxy$, and $u_nx^*$ lie in $L_m(M)$.
Thus $A$ is a $*$-algebra, and $M\subseteq A$.
Theorem~\ref{thm:sandwich} puts it inside $LS(M)$.
Every bounded operator affiliated with $M$ belongs to $M$,
so $A_b=M$. For every $x\in A$, bounded functional calculus
gives $(1+x^*x)^{-1}\in M\subseteq A$; these are precisely
the defining $EW^*$ conditions together with $A_b=M$.
In particular, by the established $EW^*$--solid correspondence,
$A$ is a solid $*$-subalgebra of $LS(M)$ containing $M$
\cite[Section~5]{BChS}.

A central element $a$ of $A$ commutes in particular with all
of $M$. The relative commutant of $M$ in $LS(M)$ is
$Z(LS(M))=S(Z(M))$: indeed, an affiliated operator commuting
with all unitaries of $M$ has central spectral projections
(and the same conclusion follows for its polar part), hence
is affiliated with $Z(M)$. Therefore $Z(A)\subseteq S(Z(M))$.
Conversely any $c\in S(Z(M))$ has bounded central spectral
cuts increasing to $1$, hence $c\in E(M)\subseteq A$ by
Theorem~\ref{thm:sandwich}; it commutes with all of
$LS(M)$. This yields \eqref{eq:center}.
\end{proof}

\begin{thm}[Innerness for general type II algebras]
\label{thm:inner}
Let $M$ be a von Neumann algebra of type II, without a
factoriality, separability, or $\sigma$-finiteness assumption,
and let $m$ be any Ciach measure. Then every additive ring
derivation
\[
 D:LL_m(M)\longrightarrow LL_m(M)
\]
is inner. In particular every such derivation is automatically
complex-linear, and
\begin{equation}\label{eq:Lie}
 Der(LL_m(M))\cong_{\mathrm{Lie}}
 LL_m(M)/S(Z(M)).
\end{equation}
The class on the right is an algebraic complex Lie quotient.
\end{thm}
\begin{proof}
Put $A=LL_m(M)$. By Lemma~\ref{lem:ew}, $A$ is an
$EW^*$-algebra with bounded part $A_b=M$. Because $M$ is of
type II, it has no finite type I direct summand. Applying
\cite[Theorem~1.2]{HKS} to the given ring derivation, we get
an implementer $a\in LS(M)$ such that
\begin{equation}\label{eq:spatial}
 D(x)=[a,x]\qquad(x\in A).
\end{equation}
At this point we only know that $a$ lies in $LS(M)$;
we must still show that an implementer can be chosen in $A$.

As established in Lemma~\ref{lem:ew}, $A$ is a solid
$*$-subalgebra of $LS(M)$.
Write $a=b+ic$ with $b,c\in LS(M)$ self-adjoint.
For every unitary $u\in M$, $[a,u]\in A$ by
\eqref{eq:spatial}. Also
$[a^*,u]=-[a,u^*]^*\in A$. Consequently
$[b,u]$ and $[c,u]$ belong to $A$.
Apply the estimate \eqref{eq:estimate} to $b$ with
$\varepsilon=1/2$. For some self-adjoint
$d_b\in S(Z(M))$ and a unitary $u_b\in M$,
$|b-d_b|\le2|[b,u_b]|$.
Since $[b,u_b]\in A$ and $A$ is solid, $b-d_b\in A$.
By \eqref{eq:center}, $d_b\in A$, so $b\in A$.
Exactly the same argument applied to $c$ gives $c\in A$.
Thus $a=b+ic\in A$, and \eqref{eq:spatial} shows that $D$
is an inner derivation of $A$ (not merely spatial in
$LS(M)$).
Since commutators are complex-linear, so is $D$.
Finally the Lie homomorphism
$A\to Der(A)$, $a\mapsto ad_a$, is surjective with kernel
$Z(A)=S(Z(M))$ by Lemma~\ref{lem:ew}. Its induced Lie
isomorphism is \eqref{eq:Lie}.
\end{proof}

\begin{rem}[Scope and citation discipline]
The factor and countable-product classifications are proved
above from Ciach's spectral criterion and projection comparison.
The arbitrary-center innerness theorem is a deduction from
the published \emph{ring derivation theorem}, the
$EW^*$-solid correspondence, and the published commutator
estimate. These three external statements are not claimed to
be newly proved here. In particular, spatial implementation
in $LS(M)$ and inner implementation in $LL_m(M)$ are
\emph{different assertions}; the second is established by
the additional solid-algebra argument.
\end{rem}

%\section{Consequences, limitations, and research perspective}

\begin{cor}[Every possible proper inclusion occurs]
\label{cor:all}
For type II von Neumann algebras and Ciach measures, all of
the following cases are realized:
$$
LL_m(M)=E(M)\subsetneq LS(M), \quad
 E(M)\subsetneq LL_m(M)\subsetneq LS(M),
$$
$$
LL_m(M)=LS(M).
$$
\end{cor}
\begin{proof}
For the first case use $m_0$ on a type $\mathrm{II}_1$
factor (or on any algebra with a nonzero type II summand),
using Theorem~\ref{thm:sandwich} and
Proposition~\ref{prop:strictE}. The second case is
Theorem~\ref{thm:two}. For the third case take a type
$\mathrm{II}_1$ factor $R$ and $m=\tau$; then
$LL_m(R)=LS(R)=S(R)$ by Theorem~\ref{thm:factor}.
\end{proof}

\begin{rem}[What the factor theorems do not imply]
The global quantity $\inf_{p\neq0}m(p)$ alone is not enough
for a nonfactor, and Theorem~\ref{thm:uncountable} shows that
even all individual atomic-corner thresholds are not enough
for arbitrary \emph{uncountable} products. Thus the countable
atomic-center theorem is not an unqualified classification of
all non-atomic centers or all direct integrals. The general
unconditional statements for type II $M$ are the universal
sandwich, $LL_{m_0}(M)=E(M)$, and the derivation theorem
\ref{thm:inner}. A finer, measurable-field-based description
for diffuse centers requires additional structure and is not
asserted here.
\end{rem}

\begin{rem}[On scientific novelty]
The main external tools---Ciach's $m$-measure theory,
central extensions, the innerness of derivations on
Murray--von Neumann algebras and on $EW^*$-algebras, and
commutator estimates---already belong to the literature.
The direct-product formulas, the mixed two-corner illustration,
and the uncountable-center obstruction are explicit deductions
and constructions in the present framework. Their research
priority would require a dedicated literature search before a
claim of originality or publication novelty is made.
\end{rem}

\end{document}